\documentclass[12pt]{article}

\usepackage{latexsym,amssymb}
\usepackage[T1]{fontenc}
\usepackage[dvips]{graphicx}
\newcommand{\z}{\mathbb Z}

\newtheorem{lem}{Lemma}[section]
\newtheorem{defn}[lem]{Definition}
\newtheorem{ex}[lem]{Example}
\newtheorem{co}[lem]{Corollary}
\newtheorem{thm}[lem]{Theorem}
\newtheorem{prop}[lem]{Proposition}

\newenvironment{proof}{\textbf{Proof.}}{\newline\hspace*{\fill}{$\Box$}}
\begin{document}
\title{Large mapping tori of free group endomorphisms}
\author{J.\,O.\,Button\\
Selwyn College\\
University of Cambridge\\
Cambridge CB3 9DQ\\
U.K.\\
\texttt{jb128@dpmms.cam.ac.uk}}
\date{}
\maketitle

\begin{abstract}
We present an algorithm which, given any finite presentation of a group as
input, will terminate with answer yes if and only if the group is large.
We use this to prove that a mapping torus of a finitely
generated free
group automorphism is large if it contains a $\z\times\z$ subgroup of infinite
index. We then extend this result to mapping tori of finitely
generated free group endomorphisms, as well as showing that such a group is
large if it contains a Baumslag-Solitar group of infinite index and has a
finite index subgroup with first Betti number at least 2. We also show that if
a group possesses a deficiency 1 presentation where one of the relators is a
commutator then it is $\z\times\z$, large or is as far as possible from being
residually finite.
\end{abstract}
%{\it Key Words and Phrases.} Mapping torus, 
%BNS invariant, Alexander polynomial.\\
%\hfill\\
%2000 {\it Mathematics Subject Classification.} Primary 57M05; 
%Secondary 20F65, 57N10.\\

\section{Introduction}

Recall that a finitely generated group $G$ is large if it has a finite index
subgroup possessing a homomorphism onto a non-abelian free group. 
The advantage of this
notion, as explained in \cite{p1} and \cite{ep2}, is that on defining
a group theoretic property to be a large property if it is preserved by
preimages, finite index subgroups and supergroups, and quotienting out by
finite normal subgroups, then $G$ being large implies that $G$ possesses
any large property $\mathcal P$, provided only that there is some finitely
generated group with $\mathcal P$. Examples of such properties which are
relevant to us 
include $G$ containing a non-abelian free subgroup, being SQ-universal (every
countable group is a subgroup of a quotient of $G$), having finite index
subgroups with arbitrarily large first Betti number, having exponential
word growth, and having subgroup growth of strict type $n^n$, which is
the largest possible growth for finitely generated groups (see \cite{luse}
Section 1.11 for definitions and proof that this is a large property). Thus
on proving that $G$ is large we obtain all these other properties for free.

There have been a range of results that give criteria for finitely
generated or finitely presented groups to be large. Starting
with B.\,Baumslag and S.\,J.\,Pride \cite{bp} which showed that groups with a 
presentation of deficiency at least 2 are large, we then have in \cite{gr}
a condition that implies this result, as well as a proof that a group with
a deficiency 1 presentation in which one of the relators is a proper power
is large. This latter result was also independently derived by St\"ohr in
\cite{st} and was followed by conditions for a group with a 
deficiency 0 presentation where some of the relators are proper powers to be
large, due to Edjvet in \cite{edj}. Then further conditions for a finitely
presented group to be large, all of which imply the Baumslag-Pride result,
are by Howie in \cite{how}, the New York Group Theory Cooperative in 
\cite{ny} Chapter IV Theorem 7 and a characterisation by Lackenby in 
\cite{lac}. However
in Section 2 we produce an algorithm which takes as input any finite
presentation and we prove that it will terminate with answer yes if and only
if the group $G$ given by that presentation is large. If $G$ is not large
then the algorithm will not terminate in general unless $G$ is finite. The
algorithm operates by converting one of the relevant criteria into a statement
about the Alexander polynomial of $G$. This can be checked by direct
calculation and then we show that $G$ being large means that this statement
is true for a finite index subgroup of $G$.

We do not address running times and only briefly consider the practicalities of
implementing the algorithm, although its simplicity should mean that it is
of use. Instead we establish the effectiveness of the algorithm by using it
to prove that a substantial class of finitely presented groups not previously
known to be large has this property. Our focus in this paper is
on groups of deficiency 1 as here the algorithm takes a suitably nice form.
Of course unlike groups of deficiency 2 or higher, not all groups of deficiency
1 are large: think of $\z$ or the soluble Baumslag-Solitar groups given by
the presentations $\langle x,y|xyx^{-1}=y^m\rangle$ for $m\in\z\backslash
\{0\}$. Other examples of non-large deficiency 1 groups were given by Pride, 
Edjvet and Howie in \cite{ep2} consisting of those Baumslag-Solitar groups
$\langle x,y|xy^lx^{-1}=y^m\rangle$ for $l,m\neq 0$ where $l$ and $m$ are
coprime, as well some HNN extensions of these, and one can find the odd further
example in the literature.

As for large groups of deficiency 1, we have already mentioned those with a
relator that is a proper power and we again have examples in \cite{ep2} with
Theorem 6 stating that the group $\langle x,y|x^ny^lx^{-n}=y^m\rangle$ for
$l,m,n\neq 0$ is large if $|n|>1$ or if $l$ and $m$ are not coprime. Further
results of a more technical nature which give largeness for some other
2-generator 1-relator presentations are in \cite{ephd} from 1984. Here it is 
asked if those groups
which are an extension of a finitely generated non-abelian free group by
$\z$ are large. They are certainly torsion free groups with a natural
deficiency 1 presentation and are also
referred to as mapping tori of finitely generated non-abelian free group
automorphisms. Some ad hoc examples of this type were shown there to be large.
More recently a naturally occurring class of deficiency 1 groups was shown to 
be large in \cite{clr}, \cite{but}, \cite{lac2}, namely fundamental groups of
compact orientable irreducible 3-manifolds  with (non-empty) boundary
consisting solely of tori (with $\langle x,y|xyx^{-1}=y^m\rangle$ for
$m=0,\pm 1$ making up the few small exceptions).

At this point it seems difficult to say convincingly either way whether groups
of deficiency 1 are generally large. In this paper we hope to offer substantial
evidence that largeness is a natural property to expect in a deficiency 1
group. Although we will display a few new groups of deficiency 1 which are not
large in Section 5, our main results are on establishing families of deficiency
1 groups which are all large. In Section 3 we prove that a mapping torus $G$
of a finitely generated non-abelian
free group automorphism is large if it contains a
$\z\times\z$ subgroup. By \cite{bf}, \cite{bfad} and \cite{brink},
these are exactly the mapping tori of finitely generated non-abelian free group
automorphisms which are not word-hyperbolic, thus this question is now
reduced to one only involving word-hyperbolic groups. However we make no use
here of the notion of word-hyperbolicity: it is the $\z\times\z$ subgroup
itself which is crucial to the proof of largeness. On combination with other
results the method also deduces that if $G$ is finitely generated but is
$F$-by-$\z$ for $F$ an infinitely generated free group then $G$ is large.

Mapping tori of finitely generated free group automorphisms appear
to make up a sizeable class of deficiency 1 groups but we can expand this
class considerably
by allowing arbitrary endomorphisms in place of automorphisms. Such groups have
been the attention of much recent research where significant progress has been
made. In particular these groups have been shown to be coherent (every finitely
generated subgroup is finitely presented) in \cite{fh}, Hopfian in \cite{gmsw}
and even residually finite in \cite{bs}. If largeness were added to 
this list (on removing the obvious small exceptions) 
then it would show that such a mapping
torus, indeed even a group which is virtually such a mapping torus, has all the
nice properties that one could reasonably hope for (but not more than this,
for instance $G$ is never LERF if the endomorphism is injective but not 
surjective so this seems like much too strong a property to expect in an
arbitrary deficiency 1 group).

In Section 4 we indicate how our proof of largeness for automorphisms 
generalises to mapping tori of finitely generated free group 
endomorphisms with a $\z\times\z$ subgroup, unless $G$ is
isomorphic to $\langle x,y|xyx^{-1}=y^{\pm 1}\rangle$. 
Unlike the case for automorphisms, this does not cover all such $G$ which are
not word-hyperbolic because on allowing endomorphisms we can now have $G$
containing other soluble Baumslag-Solitar groups. It is conjectured in 
\cite{kp} that a mapping torus of a finitely generated free group
endomorphism is word-hyperbolic if and only if it does not contain
Baumslag-Solitar subgroups. We also obtain in Section 4 largeness for 
mapping tori $G$ of finitely generated free group endomorphisms which contain a
Baumslag-Solitar subgroup provided $G$ has a finite index subgroup $H$ 
($\neq\z\times\z$) with $\beta_1(H)\geq 2$. Of course if $\beta_1(H)=1$ for
all $H$ then we would have an example of such a $G$ which is not large.
However we know of no examples apart from the soluble Baumslag-Solitar groups
themselves, and it seems believable that no other $G$ has this property, in
which case on assumption of this and the conjecture above we have largeness
for all non word-hyperbolic mapping tori of finitely generated free group
endomorphisms apart from the soluble Baumslag-Solitar groups above. 
But once again the
actual concept of word-hyperbolicity is not used anywhere.

The proof for largeness of mapping tori of finitely generated free group
endomorphisms containing $\z\times\z$ uses the fact that they are
residually finite. However the condition actually required in the proof is much
weaker and this is explored in Section 5. We introduce the concept of a
residually useless group and this has a number of equivalent definitions,
one of which is that it is finitely generated and non-abelian but has no
non-abelian finite quotients. It is merely the property of not being
residually useless which allows us to finish off the proof of largeness
of our mapping tori. Moreover introducing this definition provides us with
two advantages, one of which is theoretical and one of which is practical:
first we get to weaken significantly our hypotheses for largeness and second
it should be much quicker on being given a particular presentation to
determine that it is not residually useless by finding one non-abelian
finite quotient rather than having to prove it is residually finite.
We obtain Theorem 5.5 which states that if $G$ has a deficiency 1 presentation
in which one of its relators is a commutator then $G=\z\times\z$ or $G$ 
is residually useless with abelianisation $\z\times\z$ or
$G$ is large. In particular, on excluding the obvious small group,
the only way a group $G$ with such a
presentation can fail to be large is when it is as far away from being
residually finite as $G$ possibly can. This can ultimately be regarded as
our main result in the sense that our proof of largeness for mapping tori
involves a series of steps showing that they have finite index subgroups with
a presentation of this form. The theorem also gives us largeness for groups
$G$ with a 2-generator 1-relator presentation where the relator is a
commutator, unless $G=\z\times\z$ (which is easily detected) or $G$ is
residually useless. It is true that 2-generator 1-relator groups which are 
residually useless exist, but if we insist that the relator is a product of 
commutators then no examples are known; indeed it was only recently that
non residually finite examples of such groups were given. Moreover if the
relator is a single commutator then no examples are known that fail even to
be residually finite, so in this case being not residually useless and hence
large seems like a good bet.  
 
\section{The algorithm}

We first need to summarise the facts we require about the Alexander polynomial
and the Fox derivatives of a finitely presented group; see \cite{cfx}.
Let $G$ be given by a finite presentation 
$\langle x_1,\ldots ,x_n|r_1,\ldots ,r_m\rangle$ and let $G'$ be the derived
(commutator) subgroup of $G$, so that the abelianisation $\overline{G}=G/G'$ is
a finitely generated abelian group $\z^b\times T$, where $T$ is the torsion and
$\z^b$ is what we call the free abelianisation $ab(G)$ of $G$. Here 
$b=\beta_1(G)$ is the first Betti number of $G$ and we must have $b\leq n$.
The Alexander polynomial $\Delta_G$ of $G$ is an element of the group ring
$\z[ab(G)]$ which we think of as Laurent polynomials in $b$ variables with
integer coefficients, but it is only specified up to units which are the 
monomials. It is
defined by the following process: given the free group $F_n$ of rank $n$ with
free basis $x_1,\ldots ,x_n$ we have the free derivations (or Fox derivatives)
$D_j:\z[F_n]\rightarrow \z[F_n]$ such that $D_j(x_i)=\delta_{ij}$. We construct
the $m\times n$ Alexander matrix $(a_{ij})$ of our presentation for $G$ by
calculating $D_j(r_i)\in\z [F_n]$ and then the entry $a_{ij}$ is the image of
$D_j(r_i)$ under the natural map from $\z[F_n]$ to $\z[ab(G)]$ via $\z[G]$.
We then define the Alexander ideal of the presentation to be the ideal of
$\z[ab(G)]$ generated by the $(n-1)\times (n-1)$ minors (namely those obtained
by deleting one column and the appropriate number of rows) and the Alexander
polynomial is the highest common factor of these minors, up to units. The
utility of $\Delta_G$ is that it is independent of the finite presentation
for $G$ (at least once a basis is chosen for $ab(G)$ and we can certainly
cope with a basis change), which can be seen because it is invariant under
Tietze transformations.

If $\beta_1(G)=0$ then $\Delta_G$ is just a non-zero integer, and 
probably not an
interesting one at that, but our focus here is when $\beta_1(G)\geq 2$, in 
which case we require the definition of the relative Alexander polynomial.
Let $f:G\rightarrow ab(G)$ be the free abelianisation map and let $\theta$
be any group homomorphism of $G$ onto a free abelian group $\z^k$, so that
$k\leq b$ and $\theta=\tilde{\theta}f$. Then the Alexander polynomial
$\Delta_{G,\theta}$ relative to $\theta$ is formed in exactly the same way
with $\theta$ in place of $f$, and we can see that evaluation of 
$\Delta_G$
under $\tilde{\theta}$ divides $\Delta_{G,\theta}$. The important case 
for us is 
the polynomial $\Delta_{G,\chi}(t)$ when $\chi:G\rightarrow\z$ is a surjective
homomorphism, as then $\Delta_{G,\chi}$ carries specific and accessible
information about $\mbox{ker }\chi$. We quote the following essential fact.
\begin{prop}
Suppose $G$ is a finitely presented group and $\chi:G\rightarrow\mathbb Z$ is a
surjective homomorphism with ker $\chi=K$. Then the relative Alexander
polynomial $\Delta_{G,\chi}\in\mathbb Z[t^{\pm 1}]$ has degree $\beta_1(K)$,
with $\Delta_{G,\chi}=0$ if and only if $\beta_1(K)$ is infinite. 
\end{prop}
\begin{proof} See \cite{lic}. Here we interpret $\beta_1(K)=\infty$
as ker $\chi$ surjects to $\z^k$ for all $k\in\mathbb N$.
\end{proof}

This immediately gives us restrictions on $\beta_1(K)$.
\begin{prop}
If $G$ is finitely presented with $\beta_1(G)\geq 2$ and $K=\mbox{ker }\chi$
for $\chi$ a surjective homomorphism of $G$ to $\z$ then the degree of
$\Delta_{G,\chi}$ is at least $\beta_1(G)-1$. However if $\beta_1(G)=1$ then
$\Delta_G(1)\neq 0$ and so $\beta_1(K)$ is finite.
\end{prop}
\begin{proof}
The first part follows because if $G/K=Q$ then $\beta_1(G)\leq \beta_1(K)
+\beta_1(Q)$. For the second part, we take a presentation such that the
generator $x_1$ appears with exponent sum zero in all the relators $r_i$.
This can always be achieved by Tietze transformations and we say that such
a presentation is in standard form with respect to $x_1$. Note that 
consequently
the images of the $x_j$ in $\overline{G}$ must have finite order for $j\geq 2$.
On forming the Alexander matrix we see that all entries $D_1(r_i)$
in the first column are zero; otherwise we put $t=f(x_1)$ equal to 1
so that $D_j(r_i)$ becomes the exponent sum matrix of the presentation. Hence 
if $\Delta_G(1)=0$ then all minors evaluated by deleting the first column are 
zero, thus any $n-1$ relators are linearly dependent when abelianised, meaning 
that $\beta_1(G)\geq 2$.
\end{proof}

The next point is the crucial fact which allows us to use the Alexander
polynomial to detect largeness.
\begin{thm}
If $G$ is a finitely presented group with a homomorphism $\chi$ onto $\z$
such that $\Delta_{G,\chi}=0$ then $G$ is large.
\end{thm}
\begin{proof}
This is obtained by examining Howie's condition for largeness in \cite{how}
Section 2. Adopting that notation, we let $K$ be the standard connected
2-complex obtained from a finite presentation of $G$ consisting of $n$
generators and $m$ relators, with $N=\mbox{ker }\chi$ and $\overline{K}$ the
2-complex which is the regular covering of $K$ corresponding to $N$ so that
$\pi_1(\overline{K})=N$. Let $F$ be a field: on following
through the proof of \cite{how} Proposition 2.1, we see that if 
$H_1(\overline{K};F)$ contains a free $F[\z]$-module of rank
at least 1 then the conclusion of the proposition holds. But this is the
hypothesis of \cite{how} Theorem 2.2 which proves that for any sufficiently
large $n$ the finite index subgroup $NG^n$ admits a homomorphism onto
the free group of rank 2.

In our case we have that $H_1(\overline{K};\z)$ is a module over 
$\z[t^{\pm 1}]$ where $t$ is the generator of Im$(\chi)$ and acts by 
conjugation on $H_1(\overline{K};\z)$. The process of forming the Alexander
matrix $A$ with respect to $\chi$ using the free differential calculus results 
in a presentation matrix $P$ for the module $H_1(\overline{K};\z)\oplus\z
[t^{\pm 1}]$ (which is an $n\times m$ matrix as $P$ is actually the transpose
of $A$), with the second term in this direct sum meaning that we can make a
basis change to $P$ so that the bottom row is zero
(see \cite{lic} page 117 for details). Thus $H_1(\overline{K};\z)$ is a 
finitely presented $\z[t^{\pm 1}]$-module. Now let us move to rational
coefficients, thus taking $F=\mathbb Q$ in Howie's result above.
We still have that $H_1(\overline{K};\mathbb Q)$ is a finitely
presented module over $\mathbb Q[t^{\pm 1}]$ but this is a principal ideal
domain, so by the structure theorem it is a direct sum of cyclic modules.
Thus the presentation matrix $P$ can be put into canonical form in which all
off-diagonal entries are 0, as well as the bottom row. The Alexander polynomial
$\Delta_{G,\chi}$ (now having rational coefficients but still of the same 
degree) is just the product $p_1\ldots p_{n-1}$ where $p_i\in\mathbb Q
[t^{\pm 1}]$ are the diagonal entries, with $p_n=0$. This is because
the elementary ideals are
invariant, but the $(n-1)\times (n-1)$ minors are all zero unless the bottom
row is crossed out. As $\Delta_{G,\chi}=0$, we must have $p_i=0$ for some $i<n$
so $H_1(\overline{K};\mathbb Q)$ does have a free $\mathbb Q[\z]$-module of
rank 1 in its decomposition.
\end{proof}
\begin{co} If $G$ is a finitely presented group possessing a homomorphism to 
$\z$ with kernel having infinite Betti number then $G$ is large.
\end{co}
\begin{proof}
This is just Proposition 2.1 and Theorem 2.3. In fact Theorem (B) of the
Digression in \cite{gr} Section 4.5 has
a similar statement with bounded cohomology.
\end{proof}

Note: 1. If $n-m\geq 2$ then there are no $(n-1)\times (n-1)$ minors, whence
we define the Alexander ideal (and consequently the Alexander polynomial) to
be zero. Thus Theorem 2.3 also includes the ``deficiency at least 2 implies 
large'' result.\\ 
2. The Corollary is most definitely not true for all finitely generated groups;
we do require a finite number of relators too. Let $G$ be the restricted
wreath product $\z \wr \z$ which has presentation
\[\langle x_i,y|x_i=y^ix_0y^{-i},x_0x_ix_0^{-1}=x_i\rangle\qquad
\mbox{for $i$ in }\z.\]
This is generated by $x=x_0$ and $y$ and is soluble (with $G'$ abelian, 
although it is infinitely generated) so is a long way from being large. But the
homomorphism $\chi$ with $\chi(x)=0,\chi(y)=1$ has as its kernel the direct
product of $\z$ copies of $\z$.

Of course the converse of Corollary 2.4 is not true in general because we 
certainly have finitely presented groups $G$ which are large but with 
$\beta_1(G)\leq 1$, so $\Delta_{G,\chi}\neq 0$ by Proposition 2.2. However
such a $G$ will have a finite index subgroup $H$, for which we write 
$H\leq_f G$, that surjects to the non-abelian free group $F_2$ of rank 2 and
certainly $\beta_1(H)\geq 2$. Thus we must ask whether Corollary 2.4 recognises
largeness when we are ``staring a free group in the face'' and indeed it does.

\begin{prop} If $H$ is a finitely presented group which has a surjective
homomorphism $\theta$ to a non-abelian free group $F_n$ of rank $n\geq 2$ then
we have homomorphisms $\chi$ from $H$ onto $\z$ with $\Delta_{H,\chi}=0$.
\end{prop}
\begin{proof}
There are homomorphisms $\chi$ onto $\z$ which factor through $F_n$; take any
one of these so that $\chi=\tilde{\chi}\theta$. Then $\theta$ sends ker $\chi$
onto ker $\tilde{\chi}$, but the free group $F_n$ has no non-trivial finitely
generated normal subgroups of infinite index, so ker $\tilde{\chi}$ is an 
infinitely generated free group with $\beta_1(\mbox{ker }\tilde{\chi})=\infty$.
Thus $\beta_1(\mbox{ker }\chi)=\infty$ and $\Delta_{H,\chi}=0$
by Proposition 2.1.
\end{proof}

If we regard the space of homomorphisms from $G$ to $\z$ as $\z^{\beta_1(G)}$ 
then we can think of Proposition 2.5 as saying that we must have at least a
whole 2 dimensional ``subspace'' of homomorphisms $\chi$ with $\Delta_{G,\chi}
=0$ for $G$ to surject to a non-abelian free group (and this is not even
sufficient; see Example 2.7) whereas Theorem 2.3 only requires a
1 dimensional subspace for $G$ to be large. Thus the idea is that Theorem 2.3
should be able to
detect largeness ``at a distance'', given that we do not need to find
a specific finite index subgroup $H$ of $G$ surjecting to $F_2$ if we just 
want to establish that $G$ is large.

\begin{thm} 
There is an algorithm which, on being given any finite presentation
as input, is guaranteed to terminate with the answer yes if and only if the 
group $G$ defined by that presentation is large (but which might not terminate
if $G$ is not large).
\end{thm}
\begin{proof}
Recall that there is an algorithm which takes as input a finite presentation
and a positive integer $n$ and which outputs all the (finitely many)
subgroups $H$ having index $n$ in the group $G$ defined by the presentation. 
This is shown in \cite{disch} and is based on the Todd-Coxeter coset 
enumeration process. The output for each $H$ is a list of generators of $H$ and
a coset table for the right regular action of $G$ on the cosets of $H$. This
allows us by the Reidemeister-Schreier rewriting process to give a finite
presentation for $H$. Thus our algorithm works by taking each $n$ in turn and
each $H$ of index $n$ and using the Fox derivatives to form the Alexander 
matrix, thus enabling us to find the Alexander polynomial $\Delta_H$
(which is the highest common factor of elements in the unique factorisation
domain $\z [t^{\pm 1}_1,\ldots ,t^{\pm 1}_k]$ so we can calculate it using
\cite{coh} Algorithm 3.2.10). We now wish to see whether there is
a homomorphism $\chi$ from $H$ onto $\z$ with $\Delta_{H,\chi}=0$ which
would prove that $H$, and hence $G$, is large by Theorem 2.3. For this we
require $\beta_1(H)\geq 2$ by Proposition 2.2, so we move on if 
$\beta_1(H)\leq 1$. Otherwise we have infinitely many $\chi$ but Lemma 2.7
below shows that given $\Delta_H$, it is a finite process to determine
whether there is $\chi$ with $\Delta_{H,\chi}=0$. If there is then we stop
with the answer $G$ is large. Otherwise we move on to the next finite index
subgroup, and if $G$ is large we will eventually find an $H$ with a surjective 
homomorphism to $F_2$ so by Proposition 2.5 we will have $\Delta_{H,\chi}=0$
and we will stop with the answer yes.

If $G$ is finite with order $N$ then when $n=N$ we are just performing 
Todd-Coxeter coset enumeration with the identity subgroup, which is guaranteed
to terminate and give this order so our algorithm
stops too and would output: ``No, $G$ is not large''. 
However if $G$ is infinite
but not large then it is clear that the algorithm will run for ever (even if
$G$ has no proper subgroups of finite index, as we still need to run the
finite index subgroups algorithm for each $n$ to confirm there are no subgroups
of index $n$). 
\end{proof}

One could arrange the algorithm so that in certain special cases it 
recognises an infinite group $G$ which is not large, for instance
it might notice that all generators commute and so $G$ is abelian, allowing
it to output the answer no. However this will not determine whether a finitely
presented group is large or not for all inputs. It might be assumed that this
problem is provably insolvable, as is the case for so many group theoretic
properties. However this is not immediate and is an interesting question. The
two standard methods used in establishing unsolvability are to show that the
property is Markov or is incompatible with free products (see \cite{ls}
Chapter IV Section 4). Now every non-large finitely presented group can be
embedded in a large finitely presented group and vice versa (for we can embed
any finitely presented group in one which has no proper finite index
subgroups by \cite{bh} Part III Chapter $\Gamma$ Proposition 7.7). 
Moreover largeness and non-largeness are both easily seen to be compatible with
free products. If we call Theorem 2.6 a ``half'' or ``partial'' algorithm in
that it is proved to recognise largeness but may not terminate when faced
with the complementary property, we remark that unsolvable properties can
have half algorithms, for instance triviality and finiteness (both Markov
properties) are recognised by Todd-Coxeter coset enumeration and a half
algorithm for word hyperbolicity (also a Markov property as word hyperbolic 
groups cannot contain 
$\z\times\z$) is given in \cite{pap}. In addition it is noted 
in \cite{htrs} that
there is a (computationally inefficient) half algorithm for a finitely
presented group to have a homomorphism onto a non-abelian free group
(which is not relevant to our half algorithm because it would not terminate
if there is no homomorphism, and if there is then we would pick up largeness
more quickly by Proposition 2.5 and Theorem 2.3) but it is not known if there
is a half algorithm for the non-existence of such a homomorphism.

\begin{ex}
\end{ex}
Suppose $G=\langle x,y|r\rangle$ is a 2-generator 1-relator group such that
the exponent sum of $x$ in $r$ is zero (which can always be arranged by
Tietze transformations). A quick way to calculate $\Delta_{G,\chi}\in\z
[x^{\pm 1}]$ for $\chi$ the homomorphism sending $x$ to 1 and $y$ to 0  is to
write $r$ in terms of $y_i=x^iyx^{-i}$ for $i\in\z$ and then abelianise, with
each appearance of $y_i$ contributing a term $x^i$. This can be done quickly
by drawing out the relation a letter at a time, where an appearance of $x$
moves us up a level (and down one for $x^{-1}$), and we record the number of
letters $y$, with their sign, that appear at each level thus giving us the
coefficients of $\Delta_{G,\chi}(x)$.

Now suppose $\beta_1(G)=2$ so that the exponent sum of $y$ in $r$ is zero too.
We see that $\Delta_{G,\chi}(x)=0$ if and only if the exponent sum of $y$ at
each level is zero. By taking any $r$ for which this holds, we obtain many
2-generator 1-relator groups which are large. Assuming that $r$ is non-trivial,
these 2-generator groups cannot surject to $F_2$ directly. If we take $r\in
F_2''$ then $\Delta_G(x,y)=0$ and so $\Delta_{G,\chi}=0$ for all $\chi$.

\begin{ex}
\end{ex}
We give a brief outline of how the algorithm allows us to show that 
a particular 
2-generator 2-relator group is large. The group $G=\langle x,y\rangle$ is the 
fundamental group of the closed hyperbolic 3-manifold v1539(5,1) as taken from
the census \cite{cen}; it is the only example where $\beta_1(G)\geq 2$. On
getting Magma to take the given presentation and calculate the abelianisation
of the low index subgroups of $G$, we see two subgroups of index 5 with
abelianisation $\z^6$. On requesting a presentation of these subgroups, we find
that the first has six generators $x^2y$, $xyx$, $yx^2$, $y^2x^{-1}$, $x^{-1}
y^2$, $y^{-1}xy^{-1}$ and seven relators. We start to form the Alexander matrix
and in doing so send various generators to the identity in order to simplify
the terms. This leads us to $\chi$ with $\Delta_{G,\chi}=0$. In fact we then
see that our subgroup surjects to the free group generated by $xyx$ and
$y^{-1}xy^{-1}$ by setting all other generators to 1.\\

We make a few points on how the algorithm might be implemented before
finishing this section with the lemma we need. Firstly it would be more
efficient just to look at one subgroup from each conjugacy class and this is
exactly how the low index subgroup algorithm operates in Magma and GAP. It was
noted in Theorem 2.6 that we can only get a positive answer with a subgroup 
$H$ where $\beta_1(H)\geq 2$. As $\beta_1(H)\geq\beta_1(G)$ for $H\leq_f G$,
this is always the case if $\beta_1(G)\geq 2$. If $\beta_1(G)\leq 1$ then
these packages can calculate the abelianisation $\overline{H}$ without needing
to rewrite in order to obtain a full presentation for $H$ because they
abelianise the
relations as they go along. Thus it is probably best to only allow a finite
presentation to be input if it has first Betti number at least 2, and to run an
initial program to find such a finite index subgroup otherwise. Of course
this program might not terminate but if the group is large then it will.

Returning to the question of determining from the full Alexander polynomial
$\Delta_G$ whether $G$ has homomorphisms $\chi$ with $\Delta_{G,\chi}=0$, we
have that $\Delta_G$ is an element of $\z [t_1^{\pm 1},\ldots ,t_b^{\pm 1}]$
where $b=\beta_1(G)\geq 2$. It is helpful here to adopt the approach of 
\cite{mcm} and \cite{dunp} where we think of $\Delta_G$ as a finite subset of
lattice points in $\z^b$, with each point weighted by a non-zero integer 
obtained by regarding each monomial that appears in $\Delta_G$ with a non-zero
coefficient as a lattice point, and the coefficient as the weight. The
ambiguity of units just means that we can shift $\Delta_G$ by unit 
translations. 
\begin{lem} Given a finitely presented group $G$ with $\beta_1(G)\geq 2$,
there is an algorithm which determines whether or nor there exists a 
homomorphism $\chi$ from $G$ onto $\z$ with $\Delta_{G,\chi}=0$.
\end{lem}
\begin{proof}
Given any homomorphism $\chi$ from $G$ onto $\z$, we can evaluate
$\Delta_G$ at $\chi$ which means that we consider the one variable polynomial
$\Delta_G(t^{k_1},\ldots ,t^{k_b})$ where $\chi(t_i)=k_i$. As this is obtained
by taking highest common factors in $\z[t_1^{\pm 1},\ldots ,t_b^{\pm 1}]$ and
then evaluating, whereas the process for $\Delta_{G,\chi}$ is the other way 
round, they are both zero or non-zero together. We picture evaluation of
$\Delta_G$ at $\chi$ in the following way: on factoring $\chi$ through
$\z^b$ we have
$\tilde\chi:\z^b\rightarrow\z$ which we extend to an affine map 
$\phi:\mathbb R^b\rightarrow\mathbb R$. Then for $x\in\mathbb R$ we know that
$\phi^{-1}(x)$ is a hyperplane and $\Delta_G$ is zero on evaluation precisely
when the following condition is satisfied: for all $m\in\z$ with $\phi^{-1}(m)
\cap\Delta_G\neq\emptyset$, we require that the sum of the weights 
corresponding to the points of $\Delta_G$ in this hyperplane $\phi^{-1}(m)$ is
zero. Let us refer to this situation as ``$\Delta_G$ cancels along parallel
hyperplanes of constant $\chi$''.

We let $C$ be the convex hull of $\Delta_G$ in $\mathbb R^b$. As $\Delta_G$ is
non-empty (or else we have largeness immediately) and finite, $C$ is a convex
polytope (sometimes called the Newton polytope of the polynomial) and we refer
to \cite{bro} for details. 
We describe our algorithm by induction on the dimension and assume it for
values less than $b$. First if we are to have 
$\Delta_G(t^{k_1},\ldots ,t^{k_b})=0$ for some $\chi$ then this holds on 
putting $t=1$ so we should check that the sum of weights over all points in
$\Delta_G$ is zero and if not we are done for all $\chi$ with answer no. If so
then as well as saving time, this allows us to assume that the dimension of $C$
is $b$ as otherwise $C$ is contained in a hyperplane $H$ and we have 
cancellation purely within $H$, giving the answer yes. Now if there is some
$\chi$ such that $\Delta_G$ cancels along parallel hyperplanes of constant
$\chi$ then there must be a proper supporting hyperplane $S$ for $C$ from this 
family, giving rise to a (proper exposed) face $S\cap C$ of $C$. This cannot 
just be a vertex because $C$ being the convex hull of $\Delta_G$ means the
vertices of $C$ are contained in $\Delta_G$ but the non-zero weight on one
vertex would have nothing else to cancel out with in $S$, thus $S\cap C$ is 
also a convex polytope of dimension at least 1.

Thus we proceed by considering each 1 dimensional face (edge) of $C$ in turn
and asking whether it can lie in a supporting hyperplane $S$ giving rise to
cancellation. We check whether the sum of weights in our edge
$F_1$ is zero. If not then
such an $S$ must intersect $C$ in a 2 dimensional face $F_2$ containing $F_1$
so we replace $F_1$ with $F_2$ and repeat. Either we reach $F_{b-1}$ which now
completely determines $S$, so we can check directly for cancellation, or we
have $F_k$ with the weights in $\Delta_G\cap F_k$ having zero sum. This does
not completely determine the possible supporting hyperplane but we deal with
this by taking the quotient vector space $Q=\mathbb R^b/U$ of dimension $b-k$,
where $U$ is the subspace of dimension $k$ which is a translation of the affine
subspace generated by $F_k$. We then use the quotient map to regard $\Delta_G$ 
as a finite subset of $Q$, with new weights obtained by summing within the 
translates of $U$. Note that $\Delta_G$ cancels along hyperplanes in
$\mathbb R^b$ parallel to a supporting hyperplane containing $F_k$ if and only
if $\Delta_G$ cancels along parallel hyperplanes in $Q$, because $U$ would lie
in such a hyperplane. But by taking the convex hull of $\Delta_G$ in
$Q$ we have reduced the dimension, so we are done by induction.
\end{proof}
 
\section{Non-hyperbolic free-by-cyclic groups are large}

The deficiency of a finite presentation is the number of generators minus the
number of relators and the deficiency $def(G)$ of a finitely presented group
$G$ is the maximum deficiency over all presentations. (It is bounded above
by $\beta_1(G)$ so is finite.) We have already seen that groups of deficiency
at least 2 are large so it seems reasonable to ask whether we can use our 
algorithm to obtain large groups with lower deficiencies; clearly groups of
deficiency 1 would seem like the right place to start. In fact this turns out
to be a very fruitful choice, both from the point of view that calculating the
Alexander polynomial of a deficiency 1 group is more efficient than for lower
deficiencies, and because of the behaviour of deficiency in finite covers. To
explain this, first note that a group of lower deficiency such as the modular
group $\z_2*\z_3$ can have a finite index subgroup with deficiency at least
2 (in this case a free group) and so be proved large in this way. Thus it makes
sense to consider the virtual deficiency $vdef(G)$ of $G$ which is the supremum
of $\{def(H):H\leq_f G\}$. It is clear that $vdef(G)\geq 2$ implies that $G$
is large, but it is a good idea to divide up the possibilities into three
distinct cases: $vdef(G)\geq 2$, $vdef(G)=1$ and $vdef(G)\leq 0$. On taking
a presentation for $G$ with $n$ generators and $m$ relators where $def(G)=n-m$,
we can use Reidemeister-Schreier rewriting to obtain a presentation for an
index $i$ subgroup $H$ of $G$ with $(n-1)i+1$ generators and $mi$ relators,
thus the deficiency of $H$ is at least $(def(G)-1)i+1$. This means that
$vdef(G)\geq 2$ is equivalent to $vdef(G)=\infty$ in which case the deficiency
of a finite index subgroup tends to infinity with the index. Moreover if
$def(G)=1$ then either $vdef(G)=1$ with $def(H)=1$ for all $H\leq_f G$ or $G$
is large with $vdef(G)=\infty$.

The saving we gain with the algorithm in Theorem 2.6 when $def(G)=1$ and 
$\beta_1(G)=b\geq 2$ is that on inputting a presentation of $G$ with $n$
generators and $n-1$ relators, it appears that we need to remove $n$ columns
and calculate $n$ minors $M_i$, where we denote by $M_i$ the minor with the
$i$th column removed, then take their highest common factor in order to
calculate $\Delta_G$. In fact we can use the proof of \cite{butmp} Theorem 3.1
which shows that for $1\leq j,k\leq n$ we have
\[M_k(1-f(x_j))=M_j(1-f(x_k))\]
where $f$ is the natural ring homomorphism from $\z[F_n]$ to $\z[ab(G)]$ via
$\z[G]$. In particular, if the image of the generator $x_j$ has finite order
in the abelianisation $\overline{G}$ then there is no point in calculating
$M_j$ as it is 0, whereas if $f(x_j)\neq 0$ we can take $x_k$ such that
$1-f(x_k)$ is coprime to $1-f(x_j)$ in $\z[ab(G)]$. Thus $1-f(x_j)$ divides
$M_j$, giving $M_j=(1-f(x_j))\delta$ with $\delta$ independent
of $j$, hence $\delta$ is the highest common factor of the minors and so is
$\Delta_G$. The same also applies if we wish to calculate $\Delta_{G,\chi}$
for some particular homomorphism $\chi$, except that now we have 
$M_j=(1-t^{n_j})\Delta_{G,\chi}(t)/(1-t)$ where $n_j=\chi(x_j)$. Moreover 
we will be able to obtain deficiency 1 presentations for all finite index 
subgroups $H$ of $G$, and even if it happens that there is an $H$ with
$def(H)\geq 2$ then although the deficiency might not be picked up using the
algorithm, we will find on calculating $\Delta_H$ that it is zero so we will
obtain largeness.

A wide and important class of deficiency 1 groups is obtained by taking a free
group $F_n$ with free basis $x_1,\ldots ,x_n$ and an automorphism $\alpha$
of $F_n$ to create the mapping torus $G$ with presentation
\[\langle x_1,\ldots ,x_n,t|tx_1t^{-1}=w_1, \ldots ,tx_nt^{-1}=w_n\rangle\]
where $w_i=\alpha(x_i)\in F_n$ and thus $w_1,\ldots ,w_n$ also forms a free
basis. Equivalently $G$ is a semidirect product $F_n\rtimes\z$ which is
the same as being (finitely generated free)-by-$\z$ and implies that $G$ is
finitely generated (free-by-$\z$).
The following facts 
about such groups are well known and are summarised here. 
\begin{prop}
Let $G$ be as above then\\
(i) Each element of $G$ has a unique expression of the form
$kt^i$ where $k\in F_n$, and we multiply by the rule $k_1t^{i_1}k_2t^{i_2}
=k_1\alpha^{i_1}(k_2)t^{i_1+i_2}$.\\
(ii) For each $j\in\mathbb N$ we have the cyclic cover 
$G_j=\langle F_n,s=t^j\rangle$ of index $j$ in $G$ with presentation
\[\langle x_1,\ldots ,x_n,s|sx_1s^{-1}=\alpha^j(x_1), \ldots ,sx_ns^{-1}=
\alpha^j(x_n)\rangle.\]
(iii) If $H\leq_f G$ then $H$ is also a mapping torus of an automorphism of the
finitely generated free group $H\cap F_n$ which has finite index in $F_n$.\\
(iv) As well as those in (ii) we have many more finite index subgroups 
(assuming $F_n$ is not the
trivial group $F_0$). For instance given $F_m\leq_f F_n$ of index $d$, $\alpha$
acts on the (finitely many) subgroups of index $d$, thus on taking $\alpha^j$
which fixes $F_m$ we have $H=\langle F_m,t^j\rangle$ which is of 
index $dj$. There are other finite index subgroups but they all contain a
subgroup of this form. By using the fact that $F_n$ is
residually finite, we immediately obtain from (i), (ii) and (iv) the well
known result that $G$ is residually finite.
\end{prop}

That $G$ has deficiency no higher than 1 can be seen in a variety of ways
and we briefly mention three: the 2-complex associated to such a presentation
is aspherical so we can take the Euler characteristic; the Alexander polynomial
$\Delta_G$ is non-zero because it must divide $\Delta_{G,\chi}$, where $\chi$
is the natural homomorphism associated with the automorphism $\alpha$ given by
$\chi(t)=1$ and $\chi(x_i)=0$; or note that Proposition 3.1 (ii) gives us
subgroups $G_j$ of arbitrarily high index but with a bounded number of
generators, so $\beta_1(G_j)$ and thus $def(G_j)$ are bounded but this would
contradict the point about $vdef(G)$ if $def(G)\geq 2$. Consequently 
$vdef(G)=1$ by Proposition 3.1 (iii).

We now present the crucial point which allows us to gain largeness from
our algorithm in many cases of groups with deficiency 1.  
\begin{thm}
If $G$ is a group with a deficiency 1 presentation $\langle x_1,\ldots ,x_n|$
$r_1,\ldots ,r_{n-1}\rangle$ where one of the relators is of the form
$x_ix_jx_i^{-1}x_j^{-1}$ then $G$ is large if the subgroup of $ab(G)$ generated
by the images of $x_i$ and $x_j$ has infinite index.
\end{thm}
\begin{proof}
Note that the form of the relevant relator forces $\beta_1(G)\geq 2$ and
gives rise to a row in the Alexander matrix where there are only two non-zero
entries; these are of the form $1-f(x_j)$ in the $i$th column and $f(x_i)-1$
in the $j$th column (we can of course assume that $i\neq j$). 
To calculate $\Delta_{G,\chi}(t)$ for some 
homomorphism $\chi$ from $G$ onto $\z$ we evaluate the
terms of the Alexander matrix using $\chi$. We then cross out
a particular column, but no rows, and work out the relevant minor.
But there exists $\chi$ from $G$ onto $\z$ with both $x_i$ and $x_j$ in
the kernel, and on evaluating our row corresponding to the special relator,
our entries become $1-t^{\chi(x_j)}$ and $t^{\chi(x_i)}-1$ which are both
zero. Thus this row of zeros means that all minors are zero, so we have
$\Delta_{G,\chi}(t)=0$ and hence largeness of $G$.
\end{proof}

\begin{co}
If $G=\langle x_1,\ldots ,x_n,t\rangle$ is a mapping torus of the free group
$F_n$ with respect to the automorphism $\alpha$ such that $\alpha$ fixes a
generator $x_j$ of $F_n$, then $G$ is large if $\beta_1(G)\geq 3$ or if
$x_j$ has finite order in homology.
\end{co}
\begin{proof}
We have a relation $tx_jt^{-1}x_j^{-1}$ in a deficiency 1 presentation for $G$
so that $\beta_1(G)\geq 2$ and Theorem 3.2 applies unless $t$ and $x_j$
generate the homology up to finite index.
\end{proof}

Of course a group $G$ with a presentation having a relator of the form 
$xyx^{-1}y^{-1}$ is not equivalent to $G$ containing a $\z\times\z$ subgroup
even in the deficiency 1 case: if such a relation is present, other relations 
could imply that $x$ say has finite order (as in ($\z_2\times\z)*\z$).
As for the other way round, if $\beta_1(G)=1$ then it cannot have such a
relation, so a $\z\times\z$ subgroup cannot always be ``promoted'' to a 
defining relation without worsening the deficiency of the presentation. What we
can do in the case of a free group automorphism mapping torus is promote any
$\z\times\z$ subgroup to a relation in a finite index subgroup $H$ and this is
how we will proceed. We are then done provided $H$ has enough homology to apply
Corollary 3.3 and we can obtain this by again moving to finite index subgroups,
first so as to gain finite order homology and then we can use this to
increase the first Betti number.

Hence let $G$ be a mapping torus of an automorphism $\alpha$ of the free group
$F_n$. We say that $\alpha$ has no periodic conjugacy classes if whenever
$\alpha^j(w)$ is conjugate in $F_n$ to $w$ for $w\in F_n$, we have $w=1$ or
$j=0$. Taken together, the results of \cite{bf}, \cite{bfad} and \cite{brink} 
show that this is equivalent to $G$ containing no subgroups isomorphic to 
$\z\times\z$ and also to $G$ being word-hyperbolic. Thus we are proving here
the largeness of the non-word-hyperbolic mapping tori of automorphisms of
$F_n$ and we leave open the question of whether we have largeness in the
word-hyperbolic case, which is perhaps surprising because it might be thought
that the well developed theory of word-hyperbolic groups would mean they
are the more tractable case. However we point out that in the proof we do not
refer to word-hyperbolicity at any point; one uses existence of the 
periodic conjugacy class directly. We show here that the first two conditions
are equal because it can be established quickly by elementary means.

\begin{lem}
The mapping torus $G$ of an automorphism $\alpha$ of the free group $F_n$
contains $\z\times\z$ if and only if $\alpha$ has a periodic conjugacy
class.
\end{lem}
\begin{proof}
If there is $w\neq 1$ and $j\neq 0$ with $\alpha^j(w)=vwv^{-1}$ then the
automorphism $\gamma_v^{-1}\alpha^j$ fixes $w$, where we use $\gamma_v$ to 
denote the inner automorphism of $F_n$ that is conjugation by $v$.
Thus the mapping torus of $F_n$
with respect to $\gamma_v^{-1}\alpha^j$ contains $\z\times\z$ because if 
conjugation by $s_0$ denotes the image of $\gamma_v^{-1}\alpha^j$ then 
$w$ and $s_0$ commute, with $w^{i_1}s_0^{i_2}=1$ implying that $i_1=i_2=0$
by Proposition 3.1 (i).
But changing $\alpha^j$ to $\gamma_v^{-1}\alpha^j$ does not alter
the cyclic cover $G_j$ as a mapping torus because it is just replacing $s=t^j$ 
with $s_0=v^{-1}s$ in the presentation given by Proposition 3.1 (ii). 
Hence $\z\times\z\leq G_j\leq G$.

Now suppose $G=\langle F_n,t\rangle$ has two elements $x=kt^i$, $y=lt^j$ for
$k,l\in F_n$ which generate $\z\times\z$. We can assume $i\neq 0$ as we
cannot have both $x$ and $y$ in $F_n$, and on taking $z=x^jy^{-i}\in F_n
\backslash\{1\}$ we have that $xzx^{-1}=z$ implies $\alpha^i(z)=k^{-1}zk$.
\end{proof}

Thus our mapping torus $G$ which contains $\z\times\z$ can be assumed to
have a non-trivial element $w\in F_n$ with $\alpha(w)=w$ by dropping down to
a finite cover and using conjugation to change the automorphism.
\begin{prop}
If $\alpha$ is an automorphism of $F_n$ having $w\in F_n
\backslash\{1\}$ with $\alpha(w)=w$ and $G=\langle F_n,t\rangle$ is the 
associated mapping torus then there is a finite index subgroup $H=\langle F_m,
t^j\rangle$ of $G$ where $F_m\leq_f F_n$ has a free basis which includes $w$.
\end{prop}
\begin{proof}
We use the classic result of Marshall Hall Jnr.\,that if $L$ is a non-trivial
finitely generated subgroup of the non-abelian free group $F_n$ then there
is a finite index subgroup $F_m$ of $F_n$ such that $L$ is a free factor of 
$F_m$. We just need to put $L=\langle w\rangle$ so that $F_m=\langle w\rangle
*C$ for some $C\leq F_n$ with $w$ a basis element for $F_m$. Now we take $j>0$
with $\alpha^j$ fixing $F_m$ as in Proposition 3.1 (iv) and this gives us
$H=\langle F_m,s=t^j\rangle$ whose natural presentation has deficiency 1 and
contains the relation $sws^{-1}=w$.
\end{proof}

We are now in the position to apply Corollary 3.3 to $H$ if $\beta_1(H)\geq 3$,
obtaining largeness. But if $\beta_1(H)\geq 2$ then we need to take further
finite index subgroups in order to gain more homology. The ``smallest''
possible abelianisation of $H$ is $\z\times\z$ in which case we need to
proceed in two steps. It is only now where we require that our free
group $F_n$ is non-abelian so that $n\geq 2$. If $n=1$ then we will have 
reached this point with $H$ itself equal to $\z\times\z$, but clearly all
subgroups are of this form and $H$ is not large.
\begin{lem}
If $G$ is the mapping torus of an automorphism $\alpha$ of the non-abelian free
group $F_n$ with abelianisation $\overline{G}=\z\times\z$ then there is a
cyclic cover $G_j$ of $G$ with its abelianisation $\overline{G_j}$ having a
surjective homomorphism to $\z\times\z\times\z_m$ for some $m\geq 2$.
\end{lem}
\begin{proof}
The automorphism $\alpha$ induces an automorphism $\overline{\alpha}$ of
$\z^n$ which is an element of $GL(n,\z)$. Take any prime $p$ and
regard $\overline{\alpha}$ as an element of the finite group $GL(n,\z_p)$ then,
on taking $j$ with $\overline{\alpha}^j=1$, the cyclic cover $G_j$ has
a homomorphism onto the mapping torus of $\z_p\times\ldots\times\z_p$ formed 
from the identity automorphism. Thus $G_j$ maps onto $\z\times\z_p\times\z_p$, 
as well as onto $\z\times\z$ because $\beta_1(G_j)\geq\beta_1(G)$.
\end{proof}

We can now get a finite index subgroup with first Betti number at least 3 by
using Reidemeister-Schreier rewriting with respect to abelian covers.

\begin{prop}
If $G=\langle x_1,\ldots , x_n,t|r_1,\dots ,r_n\rangle$ has a deficiency 1
presentation with a relator $r_1=tx_1t^{-1}x_1^{-1}$ and the abelianisation
$\overline{G}=\z\times\z\times\z_m$ for $m\geq 2$ then $G$ is large.
\end{prop}
\begin{proof}
We are done by Theorem 3.2 unless the images $\overline{x_1},\overline{t}$ in
$\overline{G}$ 
generate a finite index subgroup $S$ of $\overline{G}$, but if so then $S$ can 
only be isomorphic to $\z\times\z$. Take a homomorphism $\theta$ from $G$ onto
$\z_j$ for $j\geq 2$ such that $S$ is in the kernel. We require another
generator $g\in\{x_2,\ldots ,x_n\}$ such that $\theta(g)$ generates 
Im $\theta$ but this can be achieved by Tietze transformations amongst 
$x_2,\ldots ,x_n$. We now perform 
Reidemeister-Schreier rewriting to obtain from our original presentation of $G$
a deficiency 1 presentation for ker $\theta$ consisting of $nj+1$ generators
and $nj$ relators. We have $g^i$, $0\leq i<j$ as a Schreier transversal 
for ker $\theta$ in $G$ and on setting $t_i=g^itg^{-i}$ and
$x_{1,i}=g^ix_1g^{-i}$, which will all be amongst the
generators for our presentation of ker $\theta$ given by this process (because
$t,x_1\in S\leq\mbox{ker }\theta$), our original relator $r_1$ gives rise to  
$j$ relators $t_ix_{1,i}t_i^{-1}x_{1,i}^{-1}$ in the presentation
for our subgroup. As these disappear when we abelianise, we see that 
$\beta_1(\mbox{ker }\theta)$ is at least $j+1$ and we are done by Theorem 3.2.
\end{proof}

We have now covered all cases of finitely generated non-word-hyperbolic 
free-by-cyclic groups.
\begin{thm}
If $G$ is a finitely generated group which is $F$-by-$\z$ for $F$ free 
then $G$ is large if $F$ is infinitely generated or if $G$ contains 
$\z\times\z$, with the sole exception of $F=\z$ and $G=\langle x,y|xyx^{-1}=
y^{\pm 1}\rangle$.
\end{thm}
\begin{proof}
We have by \cite{fh} that if $G$ is finitely generated then it is finitely 
presented, even if $F$ is infinitely generated. If this is so then $G$ is large
by Corollary 2.4. If $F$ is finitely generated then $G$ is a mapping torus of 
an automorphism $\alpha$ of $F_n$ which contains 
$\z\times\z$ and on applying Lemma 3.4 
and Proposition 3.5 whilst dropping to the appropriate finite index subgroups,
we can assume that $G$ has a relator of the form $tx_1t^{-1}x_1^{-1}$ in a
deficiency 1 presentation and thus $\beta_1(G)\geq 2$. If $\beta_1\geq 3$ then
apply Corollary 3.3 directly, whereas if $\beta_1(G)=2$ then use Lemma 3.6 if
necessary, noting that the cyclic cover will preserve the form of the relator,
and then Proposition 3.7.
\end{proof}

One problem that will be encountered if trying to extend this result to
word-hyperbolic groups of the form $F_n$-by-$\z$ is that it would imply
they have finite index subgroups with first Betti number at least two. 
This is unknown as we have a question of Casson which appears in various
problem lists: ``Does every automorphism $\alpha:F_n\rightarrow F_n$ leave
invariant a finite index subgroup $K$ such that $\alpha_{ab}:K/K'\rightarrow
K/K'$ has an eigenvalue which is a root of unity'' is asked in \cite{bpr}
Question 12.16 whereas in \cite{nypb} (F33) 
all eigenvalues being roots of unity are required. Note that the weaker
version is equivalent to our question because if $G$ is $F_n$-by-$\z$ with
$L\leq_f G$ having $\beta_1(L)\geq 2$ then we can drop down to a subgroup
$H\leq_f L$ of the form in Proposition 3.1 (iv), 
with $H=\langle F_m,t^j\rangle$ where $F_m\leq_f F_n$. 
Now take $K\leq_f F_m$ which is characteristic in $F_n$, then
$\alpha$ leaves $K$ invariant with $\langle K,t^j\rangle\leq_f L$ having
first Betti number at least 2, hence $\alpha_{ab}^j$ restricted to $K/K'$
has 1 as an eigenvalue. At least we see Casson's question is true in the
non-word-hyperbolic case.

However we can say something for all free-by-cyclic
groups if we relax our notion of large somewhat. Recall that a finitely
generated group $G$ is said to satisfy the Tits alternative if it contains
a non-abelian free group or is virtually soluble, although some stern
individuals insist that for a proper Tits alternative every finitely
generated subgroup of $G$ must fall into this dichotomy.
\begin{co}
If $G$ is a finitely generated group that is free-by-cyclic then either
$G$ is SQ-universal or $G$ is virtually abelian. The same is true of any
finitely generated subgroup $H$.
\end{co}
\begin{proof}
If $G$ is $F$-by-cyclic for $F$ free then for any $H\leq G$ we have that
$H$ is $(F\cap H)$-by-cyclic. If the cyclic quotient of $G$ is finite then
$G$ and $H$ are virtually free. If it is $\z$ then $G$ is large by Theorem 3.8
(thus implying SQ-universality), or $G$ is virtually $\z\times\z$ or $G$ is
word-hyperbolic. But by \cite{ol} all word-hyperbolic groups are SQ-universal
or virtually cyclic.
As for $H$, either it is contained in $F$ or it too is free-by-$\z$.   
\end{proof}

\section{Ascending HNN extensions of free groups}

Having shown that mapping tori of finitely generated non-abelian free group
automorphisms are large if they contain $\z\times\z$, we now turn to
arbitrary endomorphisms. Given a homomorphism $\theta:F_n\rightarrow F_n$ 
we can again form the mapping torus
\[G=\langle x_1,\ldots ,x_n,t|tx_1t^{-1}=\theta(x_1),\ldots ,
tx_nt^{-1}=\theta(x_n)\rangle\]
but we are not now assuming that $\theta$ is injective or surjective. However
there is a neat way of sidestepping the non-injective case using \cite{kpnt}
where it is noted that $G$ is isomorphic to a mapping torus of an injective
free group homomorphism $\tilde{\theta}:F_m\rightarrow F_m$ where $m\leq n$.
Of course it might be that $F_n$ is non-abelian but $m=0$ or 1 in which case
$G=\z$ or $\langle a,t|tat^{-1}=a^k\rangle$ for $k\neq 0$. However in these
cases $G$ is soluble and so is definitely not large. Therefore we will assume
throughout that $\theta$ is injective, in which case $G$ is also called an
ascending HNN extension of the free group $F_n$, where we conjugate the base 
$F_n$ to an isomorphic subgroup of itself using the stable
letter $t$.

We also assume in this section that $\theta$ is not surjective, whereupon we
call the ascending HNN extension proper, as automorphisms have already been
covered in the last section. We now develop the basic facts of ascending HNN
extensions of finitely generated free groups which we will find
mostly mirror the case of automorphisms, but which require more care with the 
proofs. We have that our base $F_n=\langle x_1,\ldots ,x_n\rangle$ embeds in
$G$ and we will refer to this copy of $F_n$ in $G$ as $\Gamma$, with $\theta
(\Gamma)<\Gamma$ being isomorphic to $F_n$ which means 
it has infinite index in $\Gamma$.

Once an ascending HNN extension $G$ is formed, there is an obvious homomorphism
$\chi$ from $G$ onto $\z$ associated with it which is given by $\chi(t)=1$
and $\chi(\Gamma)=0$. In the case of an automorphism ker $\chi$ is simply
$\Gamma$ but otherwise we have to consider ``backwards conjugates'', with
\[\mbox{ker }\chi=\bigcup_{i=0}^\infty t^{-i}\Gamma t^i.\]
Note $\theta(\Gamma)=t\Gamma t^{-1}<\Gamma$ so that ker $\chi$ is a strictly
ascending union of free groups, thus is infinitely generated and locally free,
but never free because $\beta_1(\mbox{ker }\chi)\leq\beta_1(\Gamma)$. One
consequence which makes this situation rather different from the case of
automorphisms is that if we are given $G$ which we are told is a mapping torus
with respect to an automorphism $\alpha$ of some group $\Gamma$ and we are
given the associated homomorphism $\chi$ then we can easily recover $\Gamma$
as it is ker $\chi$. But here we can replace $\Gamma$ with $t^k\Gamma
t^{-k}$ for any $k\in\z$ and still have a decomposition of $G$ as a proper
ascending HNN extension with the same associated homomorphism; indeed this
change does not alter the presentation which can be seen by adding new
generators $y_i=t^kx_it^{-k}$ in the presentation and then using them to
replace each $x_i$. In fact it can be seen (say by Corollary 4.2) that
any finitely generated subgroup $\Delta$ of
$\Gamma$ which contains $\theta(\Gamma)$ can replace $\Gamma$ whilst keeping  
$G$ as a proper ascending HNN extension with the same associated homomorphism
(although now the presentation for $G$ will of course change) so the base is
not even defined up to isomorphism. Whilst this means that we must be careful,
it also allows us to replace $\Gamma$ by a more convenient base in the proofs
if required.

The following result, which is Lemma 3.1 in \cite{gmsw}, allows us to recognise
ascending HNN extensions ``internally''.
\begin{lem}
A group $G$ with subgroup $\Gamma$ is an ascending HNN extension with $\Gamma$
as base if and only if there exists $t\in G$ with\\
(1) $G=\langle \Gamma,t\rangle$;\\
(2) $t^k\not\in\Gamma$ for any $k\neq 0$;\\
(3) $t\Gamma t^{-1}\leq\Gamma$.
\end{lem}
Given such an ascending HNN extension $G=\langle \Gamma,t\rangle$ where the
base $\Gamma$ is finitely generated, the next lemma reveals all the possible 
finitely generated bases
for expressing $G$ as an ascending HNN extension without changing the
associated homomorphism.
\begin{lem} The finitely generated subgroup $\Delta$ of $G=\langle\Gamma,t
\rangle$ is a base for the decomposition of $G$ as an ascending HNN extension
with the same associated homomorphism $\chi$ if and only if $t\Delta t^{-1}$
is contained in $\Delta$ and there are integers $k\geq l$ with
\[t^k\Gamma t^{-k}\leq\Delta\leq t^l\Gamma t^{-l}.\]
\end{lem}
\begin{proof} If so then conditions (1)--(3) in Lemma 4.1 apply to $\Delta$,
whereas if $\Delta$ is a base then $t\Delta t^{-1}\leq\Delta$ is required,
and the associated homomorphism $\chi$ would be such that
\[ \mbox{ker }\chi=\bigcup_{i=0}^\infty t^{-i}\Delta t^i\]
but as $\Gamma\leq\mbox{ker }\chi$ is finitely generated we must have $k$ with
$\Gamma\leq t^{-k}\Delta t^k$. Now swap $\Gamma$ and $\Delta$.
\end{proof}

We can now offer the equivalent of Proposition 3.1 for injective endomorphisms.
\begin{prop}
Let $G$ be a proper ascending HNN extension
\[\langle x_1,\ldots ,x_n,t|tx_1t^{-1}=\theta(x_1),\ldots ,
tx_nt^{-1}=\theta(x_n)\rangle\]
with respect to the injective endomorphism $\theta$ of the finitely generated
free group $\Gamma=F_n$ with free basis $x_1,\ldots ,x_n$
and let $\chi$ be the associated homomorphism.\\
(i) Each element $g$ of $G$ has an expression of the form 
$g=t^{-p}\gamma t^q$ for $p,q\geq 0$.\\
(ii) For each $j\in\mathbb N$ we have the cyclic cover 
$G_j=\langle \Gamma,s=t^j\rangle$ of index $j$ in $G$ with presentation
\[\langle x_1,\ldots ,x_n,s|sx_1s^{-1}=\theta^j(x_1), \ldots ,sx_ns^{-1}=
\theta^j(x_n)\rangle.\]
(iii) If $H\leq_f G$ then $H$ is also a proper ascending HNN extension of a
finitely generated free group $\Delta\leq_f\Gamma$ with respect to the
(restriction to $H$ of the) same associated homomorphism $\chi$.\\
(iv) If $\Delta\leq_f\Gamma$ then $H=\langle\Delta,t\rangle$ has finite index 
in $G=\langle\Gamma,t\rangle$.
\end{prop}
\begin{proof}
(i) is \cite{fh} Lemma 2.2 (1). The point is that positive powers can be moved
to the right and negative powers to the left, thus allowing us to multiply two
elements together. However the expression is not unique as $g=t^{-p}\gamma t^q=
t^{-p-1}\theta(\gamma)t^{q+1}$. Note that $\chi(g)=q-p$.\\
(ii) This is also Lemma 2.2 (1), but in \cite{kp}.\\
(iii) We try following the standard proof for automorphisms, but with more 
care. Setting $K=\mbox{ker }\chi$ we have $H/(H\cap K)=\z$ so let $x\in H$ be a
generator of $\z$, which we can take to be of the form $t^{-p}\gamma t^{p+m}$
for $\gamma\in\Gamma$, $p\geq 0$ and $m>0$ (where $m$ is the minimum positive
value of $\chi$ when restricted to $H$). We now put $t^{-p}\Gamma t^p$ 
in place of $\Gamma$ so we can write $x$ as $\gamma t^m$. Then setting 
$\Delta=H\cap\Gamma$ we have that $\chi(x)=m$ and $\chi(\Delta)=0$ so we 
are done if the three conditions of Lemma 4.1 are satisfied for $\Delta$
and $x$. We certainly have that $x^k\not\in\Delta$ for 
$k\neq 0$. Also $x\Delta x^{-1}\leq\Delta$ because $x\in H$ and 
$x\Gamma x^{-1}\leq\Gamma$.

In order to get that $x$ and $\Delta$ generate $H$, we show by induction that
$x^{-i}\Gamma x^i$, which is contained in $t^{-im}\Gamma t^{im}$, also
contains it. Assuming this to be true for $i$, we have
\[x^{-1}t^{-im}\Gamma t^{im}x=t^{-(i+1)m}(\theta^{im}(\gamma^{-1})
\Gamma\theta^{im}(\gamma))t^{(i+1)m}=t^{-(i+1)m}\Gamma t^{(i+1)m}.\]
Now we can write $K$ as a smaller ascending union
\[\bigcup_{i=0}^\infty t^{-im}\Gamma t^{im}=
  \bigcup_{i=0}^\infty x^{-i}\Gamma x^i\]
so that $H\cap K$ is the union of the $H\cap x^{-i}\Gamma x^i$ which equals
$x^{-i}(H\cap\Gamma)x^i$ because $x\in H$. Thus $H$ is generated by $x$ and
$\Delta$. Finally we change $\Gamma$ back to $t^p\Gamma t^{-p}$ so
that $\Delta$ is changed to $t^p\Delta t^{-p}$ which is in the original 
$\Gamma$ and which is also a perfectly good base for $H$ by Lemma 4.2.\\
(iv) Let $\gamma_1,\ldots ,\gamma_d$ be a transversal for $\Delta$ in $\Gamma$.
If we can show that $H\cap K\leq_f K$ for $K$ the kernel of the associated
homomorphism then we are done, despite the fact that
$K$ and $H\cap K$ are infinitely generated, because $t\in H$ and any $g\in G$
is of the form $kt^m$ for $k\in K$.

The set
\[S=\{t^{-m}\gamma_it^m:m\in\mathbb N,1\leq i\leq d\}\]
contains an element of every coset of $H\cap K$ in $K$. This can be seen by
writing $k\in K$ as $t^{-m}\gamma t^m$ for $\gamma\in\Gamma$ using (i). Then
there is $\gamma_i$ such that $\gamma\gamma_i=\delta\in\Delta$. This means that
$kt^{-m}\gamma_it^m=t^{-m}\delta t^m$ which is in $H$ and in $K$. We now show
that the index of $H\cap K$ in $K$ is at most $d$. Note that for $q>p$,
any element of the form $t^{-p}\gamma_it^p$ is in the same coset as some 
element of the form $t^{-q}\gamma_jt^q$ because $\theta^{q-p}(\gamma_i)
\gamma_j^{-1}=\delta$ for some $j\in\{1,\ldots ,d\}$ and some 
$\delta\in\Delta$,
thus giving $t^{-p}\gamma_it^p(t^{-q}\gamma_jt^q)^{-1}=t^{-q}\delta t^q$ which
is in $H\cap K$. Therefore we proceed as follows: $S$ is a set indexed by
$(l,i)\in\mathbb N\times\z_d$ and we refer to $l$ as the level. Choose a 
transversal $T$ for $H\cap K$ in $K$ from $S$ which a priori could be
infinite and let $g_1$ be the element in
$T$ with smallest level $l_1$ (and smallest $i$ if necessary). Then for each
level $l_1+1,l_1+2,\ldots$ above $l_1$ there is an element in $S$ with this
level that is in the same coset of $H\cap K$ as $g_1$ and so cannot be in $T$.
Cross these elements off from $S$ and now take the next element $g_2$ in $T$
according to our ordering of $S$. Certainly $g_2$ with level $l_2$ has not been
crossed off and we repeat the process of removing one element in each level 
above $l_2$; as these are in the same coset as $g_2$ they too have not been 
erased
already. Now note that we can go no further than $g_d$ because then we will
have crossed off all elements from all levels above $l_d$; thus we must have a
transversal for $H\cap K$ in $K$ of no more than $d$ elements.
\end{proof}

We now consider periodic conjugacy classes in the injective endomorphism
case. We follow \cite{kp} in generalising the previous definition of a
periodic conjugacy class which was applied to automorphisms. Let $G=\langle
F_n,t\rangle$ be the mapping torus of an injective endomorphism $\theta$ of
the free group $F_n$. We say that $\theta$ has a periodic conjugacy class 
if there exists $i>0$, $k\in\z$ and $w\in F_n\backslash
\{1\}$ such that $\theta^i(w)$ is
conjugate to $w^k$ in $F_n$. If this is so with $\theta^i(w)=vw^kv^{-1}$ then
on taking the endomorphism $\phi$ of $F_n$ such that $\phi=\gamma_v^{-1}
\theta^i$ we have on setting $\Delta=\langle w\rangle$ and $s=v^{-1}t^i$ 
that the subgroup $\langle\Delta,s\rangle$ of $G$ is an ascending HNN
extension with base $\Delta$ and stable letter $s$ by Lemma 4.1. Consequently
it has the presentation $\langle s,w|sws^{-1}=w^k\rangle$. 
These presentations are part of the
famous family of 2-generator 1-relator subgroups known as the Baumslag-Solitar
groups. We define the Baumslag-Solitar group $B(j,k)=\langle x,y|xy^jx^{-1}=
y^k\rangle$ for $j,k\neq 0$ and without loss of generality $|k|\geq j>0$ 
because we can replace $x$ with $x^{-1}$ and put it on the other side, 
as well as taking the inverse of the relation. The Baumslag-Solitar groups are
used as counterexamples in many areas of group theory so it is useful to be
aware of when they can be subgroups of mapping tori of free groups. The next
Proposition is similar to \cite{kp} Lemma 2.3 but has a completely general
conclusion. 

\begin{prop} 
Let $G=\langle F_n,t\rangle$ be the mapping torus of an injective endomorphism 
$\theta$ of the free group $F_n$.\\
(i) $G$ cannot contain a subgroup isomorphic to a Baumslag-Solitar group
$B(j,k)$ unless $j=1$ or $j=k$.\\
(ii) If there exists $i,j>0$, $k\in\z$ and $w\in F_n\backslash\{ 1\}$ with
$\theta^i(w^j)$ conjugate in $F_n$ to $w^k$ then $k=dj$ and $\theta^i(w)$ is
conjugate to $w^d$ so that $\theta$  has a periodic conjugacy class.\\ 
(iii) $G$ has Baumslag-Solitar subgroups if and only if $\theta$ has periodic
conjugacy classes. In particular $G$ contains $B(1,d)$ if and only if $G$ has
a periodic conjugacy class of the form $w\in F_n\backslash\{1\}$ and $i>0$
with $\theta^i(w)$ conjugate to $w^d$.\\
(iv) If $\theta$ is an automorphism then $G$ can only contain Baumslag-Solitar
subgroups of the form $B(1,\pm 1)$ or $B(j,j)$. This happens if and only if
$G$ has a periodic conjugacy class of the form $w\in F_n\backslash\{1\}$ and 
$i>0$ with $\theta^i(w)$ conjugate to $w^{\pm 1}$, 
so that our current definition of periodic conjugacy classes is
equivalent to the definition given for automorphisms in Section 3.
\end{prop}
\begin{proof}
We prove (ii) first. Suppose we have $\theta^i(w^j)=vw^kv^{-1}$ for $v\in F_n$
then we can do as we did in Lemma 3.4 and replace $\theta$ with the injective
endomorphism $\gamma_v^{-1}\theta^i$, thus replacing $G$ with a finite index
subgroup
of itself by Proposition 4.3 (ii). So we now have $\theta(w)=u\in F_n\backslash
\{1\}$ and $\theta(w^j)=w^k$. We set $w=c^m$ where $c$ generates a maximal
cyclic subgroup. Then $(\theta(c^j))^m=(c^k)^m$ but we are in a free group,
so $\theta(c^j)=c^k$. Hence $\theta(c)=c^d$ must also be a power of $c$ as
its $j$th power is. Thus $k=dj$ and $\theta(w)=w^d$.

For (i), if $B(j,k)\leq G$ then we have elements $x=t^{-p}at^q$, 
$y=t^{-r}bt^s$, with $p,q,r,s\geq 0$ and $a,b\in F_n$, which satisfy 
$xy^jx^{-1}=y^k$. We can 
assume that $q-p\geq 0$ by replacing $x$ with $x^{-1}$ and swapping $j$ and
$k$ (inverting the relation if now $j<0$). Then applying to both sides of the
relation the 
associated homomorphism $\chi$ of the ascending HNN extension, we have $\chi(y)
=0$ or $j=k$. With the former we have $y=t^{-r}bt^r$ so by replacing 
$H=\langle x,y\rangle$ with the conjugate subgroup $t^NHt^{-N}$ where $N=p+r$,
we obtain $y=t^pbt^{-p}$ and $x=t^rat^{q-p-r}=a't^{q-p}$ for $a'\in F_n$. Now
$q-p\neq 0$ or else $H\leq F_n$ would be cyclic as it has a relation between 
its two generators. Thus $q-p>0$ and $\theta^{q-p}(y^j)$ is conjugate via $a'$
to $y^k$, so by (ii) $\theta^{q-p}(y)$ is conjugate via $a'$ to $y^d$ which
implies from before that $H$ is isomorphic to $B(1,d)$. 

As for $B(j,j)=\langle x,y|xy^jx^{-1}=y^j\rangle$, note that this is a mapping
torus $G=\langle F_j,y\rangle$ of a free group automorphism $\alpha$ of $F_j=
\langle x=x_1,\ldots ,x_j\rangle$ where $\alpha(x_i)=x_{i+1}$ mod $j$. However
$B(1,1)=\z\times\z\leq B(j,j)$, generated by $x$ and $y^j$ (as they commute
and generate $\z\times\z$ in homology).

For (iii) we have already seen that $\theta$ having a periodic conjugacy class 
of the form $\theta^i(w)$ conjugate to $w^d$ gives rise to a subgroup $B(1,d)$
of $G$. If $B(j,k)\leq G$ then we have just said that if $j=k\geq 2$ we can
replace it with $B(1,1)$ by swapping $y$ for $y^j$. Thus following through the
proof of (i), we have for our $x$ and $y$ that generate $B(j,k)$ that either
$\chi(y)=0$ whereupon we saw that $B(j,k)$ is
actually $B(1,d)$ with $\theta^{q-p}(y)$ conjugate to $y^d$, or $j=k=1$. 
In this
case we set $l=q-p$ and $m=s-r$, both of which can be taken as non-negative by
replacing $x$ or $y$ by its inverse if necessary.
On moving positive and negative
powers of $t$ to the right and left respectively we have that $z=x^my^{-l}$ is
of the form $t^{-pm-sl}ct^{qm+rl}$ which as $c\in F_n$ means that $\chi(z)=0$. 
So as above we have some large $N$ such that we can replace $B(1,1)$ by 
$t^NB(1,1)t^{-N}$, thus we are able to regard $z$ as if it is in $F_n$, 
and we now have $x=a't^l$ and $y=b't^m$ for $a',b'\in F_n$. Now we can take
$l>0$ because $xz=zx$ so $l=0$ implies that $x$ and $z$ are in 
the same cyclic subgroup of $F_n$. But if so then swap $x$ with $y$, which
together generate $\z\times\z$. Hence $\theta$
has a periodic conjugacy class with $\theta^l(z)$ conjugate to $z$.

For (iv) we use the classic result of Higman which says that an
automorphism of a free group that maps a finitely generated subgroup into
itself maps it onto itself. Thus if $\theta$ has a periodic conjugacy class
with the automorphism $\gamma_v^{-1}\theta^i$ sending $w$ to $w^d$ then we must
have $d=\pm 1$. This happens if and only if $G$ contains $B(1,\pm 1)$ by (iii),
and we know $G$ containing $B(j,j)$ implies that $G$ contains $B(1,1)$. Now
if $\theta^i(w)=vw^{-1}v^{-1}$ then $\theta^{2i}(w)=uwu^{-1}$ where $u=
\theta^i(v)v$. Thus by (iii) $G$ contains $\z\times\z$ and this is equivalent 
to the original definition of possessing a periodic conjugacy class. 
\end{proof}

Hence just as in Section 3 we can assume that our mapping torus $G=\langle
F_n,t\rangle$ which contains $\z\times \z$ has $w\in F_n\backslash\{1\}$ with
$\theta(w)=w$. We would like to show that $G$ is large by similar means where
we first lift $w$ to a generator in a free group and then create extra finite
homology in order to apply Proposition 3.7. We shall see that the first task
can be achieved by careful but essentially elementary means as before, whereas
the second will require extra knowledge.

\begin{thm}
If $\theta$ is an injective endomorphism of the free group $\Gamma$ of rank
$n$ with $w\in F_n\backslash\{1\}$ such that $\theta(w)=w$ then there is a
finite index subgroup $\Delta$ of $\Gamma$ and $j\geq 1$ such that $\Delta$
has a free basis including $w$, with $\theta^j(\Delta)\leq\Delta$.
\end{thm}
\begin{proof}
The free basis is easy to obtain by taking $F\leq_f\Gamma$ with $\langle w
\rangle$ a free factor of $F$, just as in Proposition 3.5. The second condition
is the important part. The aim is to repeatedly pull back $F$; although we do
not have $F\leq\theta^{-1}(F)$ in general as this is equivalent to $\theta(F)
\leq F$ which would mean we are done, we do find that the index is 
non-increasing. To see this note that $\theta^{-1}(F)=\theta^{-1}(F\cap\theta
(\Gamma))$ and $\theta^{-1}\theta(\Gamma)=\Gamma$ as $\theta:\Gamma\rightarrow
\theta(\Gamma)$ is an isomorphism. Hence the index of 
$F\cap\theta(\Gamma)$ in $\theta(\Gamma)$ is preserved
by applying $\theta^{-1}$ to both sides, so is equal to the index of 
$\theta^{-1}(F)$ in $\Gamma$. But the index of $F\cap\theta(\Gamma)$ in 
$\theta(\Gamma)$ is no more than that of $F$ in $\Gamma$, thus
$[\Gamma:\theta^{-i}(F)]$ gives us a non-increasing sequence
which must stabilise at $N$ with value $k$. When it does we have for $i\geq 0$
that $\theta^{-(i+N)}(F)$ is just moving around the finitely many index $k$
subgroups. Although it happens that $\theta^{-1}$ is not in general a 
permutation of these index $k$ subgroups, we must land on some such subgroup
$\Delta$ twice so we have $j\geq 1$ with $\theta^{-j}(\Delta)=\Delta$, giving
$\Delta\geq\theta^j(\Delta)$.

We now show that, although the rank of $\theta^{-i}(F)$ reduces whenever the
index reduces, we can keep $w$ as an element of a free basis each time we pull
back. This time we restrict $\theta$ to an injective homomorphism from 
$\theta^{-1}(F)$ to $F$ with image $\theta\theta^{-1}(F)$. As $\theta
\theta^{-1}(F)$ is a finitely generated subgroup of $F$ containing a free basis
element $w$ of $F$, we can ensure $w$ is in a free basis for $\theta\theta^{-1}
(F)$ (for instance see \cite{ls} Proposition I.3.19). Now $\theta^{-1}(F)$ and 
$\theta\theta^{-1}(F)$ are isomorphic via $\theta$ with inverse $\phi$ say, 
so a basis $b_1,\ldots ,b_r$ for the latter gives rise to a basis $\phi(b_1),
\ldots ,\phi(b_r)$ for $\theta^{-1}(F)$ and if $b_1=w$ then $\phi(b_1)=w$.
\end{proof}
\begin{co}
If $G=\langle\Gamma,t\rangle$ is a mapping torus of an injective endomorphism
$\theta$ of the free group $\Gamma$ of rank $n$ and 
$\z\times\z\leq G$ then we have $H\leq_f G$ such that $H$ has a deficiency 1
presentation $\langle x_1,\ldots ,x_m,s|r_1,\ldots ,r_m\rangle$ including a
relator of the form $sx_1s^{-1}x_1^{-1}$.
\end{co}
\begin{proof}
By Proposition 4.4 (iii) we can on dropping to a finite index subgroup of $G$ 
assume that there is $w\in\Gamma\backslash\{1\}$ with $\theta(w)=w$ and then by
Theorem 4.5 we have $\Delta\leq_f\Gamma$ with $\Delta$ having a free basis
$w,x_2,\ldots ,x_m$ and $j\geq 1$ with $\theta^j(\Delta)\leq\Delta$. Thus by
Proposition 4.3 (ii) and (iv) we have that $L=\langle\Delta,s=t^j\rangle$ has
finite index in $G$ and by Lemma 4.1 $L$ is an ascending HNN extension with 
base $\Delta$ and stable letter $s$. Thus on taking the standard presentation 
for $L$ we see that it has deficiency 1 
with a relator equal to $sws^{-1}w^{-1}$.
\end{proof}

Thus such a $G$ has a finite index subgroup $H$ with $\beta_1(H)\geq 2$ and
Proposition 3.7 tells us that $G$ is large apart from in one circumstance;
namely when $\overline{H}=\z\times\z$. However, unlike the case of 
automorphisms where we used cyclic covers in Lemma 3.6 to get extra finite 
homology, we cannot do this for general endomorphisms with cyclic, abelian or 
even soluble covers as the next example shows.
\begin{ex}
\end{ex} 
Let $G=\langle F_2,t\rangle$ be the mapping torus of the free group $F_2=
\langle x,y\rangle$ (using capital letters for inverses) with respect to
the endomorphism
\begin{eqnarray*}
\theta(x)&=&x\\
\theta(y)&=&xyXYx^2Y^2X^2y^3xYXY^2x^2y^2X^2.
\end{eqnarray*}
The image $w(x,y)$ of $y$ is the commutator of the commutators $xyXY$ and
$x^2Y^2X^2y^2$. As $w$ is in $F_2'\backslash\{1\}$, it is clear that $\theta$
is injective but not surjective.
\begin{prop}
This group $G$ has the property that if $N\leq_f H\leq_f G$ with $N\unlhd G$
and $G/N$ is soluble then $G'\leq N$, so that $H\unlhd G$ with $G/H$ abelian.
Moreover $H/H'=\z\times\z$.
\end{prop}
\begin{proof}
We have that $\overline{G}=\z\times\z$ (generated by $t$ and $x$) 
and we calculate the
Alexander polynomial $\Delta_G(t,x)$. This turns out to be 1 because $w(x,y)$
has been chosen so that the group $\langle x,y|w\rangle$ has Alexander
polynomial 0 as $w$ is a commutator of commutators. 
By the same process we can see that the cyclic
covers $G_j=\langle F_2,s=t^j\rangle$ also have $\overline{G_j}=\z\times\z$
and $\Delta_{G_j}=1$ because whenever we apply $\theta$ we replace in 
$\theta^j(y)$ each appearance of $y$ with $w$, so $\langle x,y|\theta^j(y)
\rangle$ also has zero Alexander polynomial. We now show that for any 
``rectangular'' abelian cover 
$G_{j,k}$, namely the subgroup of index $jk$ in $G$ consisting of those
elements with exponent sum equal to 0 modulo $j$ in $t$ and 0 modulo $k$ in
$x$, we have $\overline{G_{j,k}}=\z\times\z$. Using Reidemeister-Schreier to
get from our presentation for $G_j$ to that of $G_{j,k}$, we have generators
$z=x^k$, $s_i$ and $y_i$ for $0\leq i\leq k-1$, where $s_i=x^isx^{-i}$ and
equivalently for $y_i$. The $k$ new relators obtained from $sxSX$ in the
presentation for $G_j$ collapse to $s_i=s$ and $szSZ$. As for the other relator
$r$ in $G_j$, we write it as $r=w_j(x,y)sYS$ where $w_j=\theta^j(y)$ has 0
exponent sum in $y$ for each level of $x$. Then in
$G_{j,k}$ we find that $r$ becomes, when abelianised, simply $Y_0$ and the
other relators for $G_{j,k}$ obtained from conjugates $x^irX^i$ are similarly
$Y_i$. Thus the $y_i$ are trivial in homology, giving us that 
$\overline{G_{j,k}}$ must be generated by $s$ and $z$, so equals $\z\times\z$.

On taking any $H\leq_f G$ we have a natural surjective homomorphism from the
abelianisation $H/H'$ of $H$ to $H/H\cap G'$. So if we set $H$ equal to
$G_{j,k}$ we have $G'\leq G_{j,k}$, meaning that this
surjective homomorphism is from $G_{j,k}/G_{j,k}'$ to $G_{j,k}/
G'\leq_f G/G'=\z\times\z$. As $\z\times\z$ is Hopfian, we
must have an isomorphism with $G_{j,k}'=G'$. Now if $A$ is any finite abelian
cover of $G$ then there is $G_{j,k}\leq A$, but $G_{j,k}'\leq A'\leq G'$ so
$A/A'=A/G'\leq_f G/G'$, giving $A/A'=\z\times\z$.

Finally if $N\unlhd_f G$ is a soluble cover then on taking $K$ with $N\unlhd_f
K\unlhd_f G$ and $(G/N)/(K/N)$ abelian, we have $K/K'=\z\times\z$ with $K'=G'$.
Now we can replace $G$ with $K$ and $K$ with some subgroup $L$ which is
normal in $K$ with $K/L$ abelian. This allows us to repeat the argument 
above because we have $G'=K'\leq L\leq_f G$ so actually $L\unlhd_f G$ with
$G/L$ abelian, thus again $L'=G'$ from above.
We continue until we reach $N$, so in fact the soluble cover was abelian with 
$N'=G'\leq N$ and $N/N'\leq_f G/G'=\z\times\z$.
\end{proof}

Therefore we require other covers in order to proceed. It was recently shown
in \cite{bs} that mapping tori of free group endomorphisms are residually
finite by using group schemes and other techniques in algebraic geometry.
This can be combined with a short lemma in \cite{lr} to complete our result on
largeness of mapping tori of injective endomorphisms of free groups. However
in proving this, it will be enough just to have a non-abelian finite quotient.
Although we will still use the result of \cite{bs} to ensure this happens,
our requirement is much weaker than the full strength of residual finiteness
and it is worthwhile investigating this condition in more detail. Therefore we
continue with the last part of our proof of largeness seemingly
unresolved, but this will be concluded in the next section.

We finish this section by looking at those mapping tori $G$ of endomorphisms
of free groups which contain a Baumslag-Solitar subgroup. Our result on
largeness is not quite definitive because we need $\beta_1(G)\geq 2$ in
order to apply our methods and we cannot show that $G$ necessarily has a
finite index subgroup with that property. However this is the only obstacle
to largeness.
\begin{thm}
If $G=\langle\Gamma,t\rangle$ is a mapping torus of an endomorphism
$\theta$ of the free group $\Gamma$ of rank $n$ which contains a 
Baumslag-Solitar subgroup $B(j,k)$ then either $G$ is large or $G=B(1,k)$ 
or $\beta_1(H)=1$ for all $H\leq_f G$.
\end{thm}
\begin{proof}
As usual we assume that $\theta$ is injective. By Proposition 4.4 we know that
$G$ can only contain Baumslag-Solitar subgroups of type $B(1,k)$ or $B(k,k)$
for $k\neq 0$ and as we have already covered those which contain $B(1,1)$, we
need only consider $B(1,k)\leq G$ for $k\neq\pm 1$. If there is some
$H\leq_f G$ with $\beta_1(H)\geq 2$ then we can replace $G$ by $H$ because
$H$ is a mapping torus by Proposition 4.3 (iii) and $B(1,k)\cap H
\leq_fB(1,k)$ so $H$ contains a Baumslag-Solitar group of the same form
(for instance use \cite{wil} which characterises these groups as the soluble
groups of deficiency 1). Therefore by Proposition 4.4 (iii) we are looking
at the case where we have a periodic conjugacy class of the form $w\in F_n
\backslash\{1\}$ and $i>0$ with $\theta^i(w)$ conjugate to $w^d$ for some
$d\neq\pm 1$. Just as in 
the $\z\times\z$ case, we drop down to a finite index subgroup and assume that
$\theta(w)=w^d$. Now we follow the proof of Theorem 4.5 to get $F\leq_f\Gamma$
with $\langle w\rangle$ a free factor of $F$, observing that $w\in\theta^{-1}
(F)$ so that we keep $w$ as we pull back $F$. Note that we can assume $w$ is
not a proper power by Proposition 4.4 (ii), so 
we can also preserve $w$ in a 
free basis each time because $w^d\in\theta\theta^{-1}(F)$ and if $w^c\in
\theta\theta^{-1}(F)$ for $0<|c|<|d|$ then the element $u\in\theta^{-1}(F)$
with $\theta(u)=w^c$ cannot be a power of $w$ but $\theta(u^d)=\theta(w^c)$,
thus contradicting injectiveness. Thus $w^d$ can be extended to a free basis
for $\theta\theta^{-1}(F)$ by \cite{ls} Proposition I.3.7
and thus $w$ will be in the corresponding basis
for $\theta^{-1}(F)$.

This gives an equivalent version of Corollary 4.6 where we have a deficiency
1 presentation with a relation $sxSX^d$. Thus on taking a surjective
homomorphism $\chi$ to $\z$ (which must send $x$ to 0) we have as in Theorem 
3.2 a top row with only one non-zero entry which is $t^{\chi(s)}-d\in\z
[t^{\pm 1}]$. If $\beta_1(G)=1$ then the only available $\chi$ will send $s$ to
$\pm 1$ but if not then we can find $\chi$ with $\chi(s)=0$. Now take a prime
$p$ dividing $1-d$ (and if $d$ is rude enough to be 2 then take the double
cover of $G$ with relation $sxSX^4$). We then obtain largeness from Howie's
criterion in Theorem 2.3 but with $\z_p$ as our field.
\end{proof}

Although we do not have a proof that a mapping torus of a 
free group endomorphism containing a Baumslag-Solitar subgroup of infinite
index has a finite index subgroup with first Betti number at least two, the
statement of Theorem 4.9 is still useful in a practical sense because if
we are presented with a particular group $G$ of this form that we would
like to prove is large, we can enter the presentation into a computer and
ask for the abelianisation of its low index subgroups. As soon as we see one
with first Betti number at least two, we can conclude largeness.
\begin{ex}
\end{ex}
The group $G=\langle a,t|t^2at^{-2}=a^2\rangle$ was shown in \cite{ds} to be
a mapping torus of a free group endomorphism (just put $b=tat^{-1}$) which
is 1-related and residually finite but not linear. It clearly contains
$B(1,2)$ so we can conclude by Theorem 4.9 that it is a large 1-related
residually finite non-linear group by getting Magma or GAP to tell us that it
has a subgroup of index six with abelianisation $\z\times\z$.\\

We can even say something if $G$ is a mapping torus of an injective
endomorphism of an infinitely generated free group in the case when $G$ is
finitely generated, thanks to the power of \cite{fh} Theorem 1.2 which proves
that $G$ has a presentation of the form
\[\langle t,a_1,\ldots ,a_k,b_1,\ldots ,b_l|ta_1t^{-1}=w_1,\ldots ,
ta_kt^{-1}=w_k\rangle\]
for $w_1,\ldots ,w_k$ words in the $a_i$ and the $b_j$. Thus either $G$ has
deficiency at least two and so is large, or $l=0$ in which case $G$ is also
a mapping torus of a finitely generated free group endomorphism and so the
results of this section apply. 

\section{Residually useless groups}
Recall that a group $G$ is residually finite if the intersection $S$ over all
the finite index subgroups $F\leq_f G$ is the trivial group $I$. Although this
works perfectly well as a general definition, it is most useful when $G$ is
finitely generated and that will be our assumption in this section. Our
motivation for the next definition is to ask: how badly can a group fail to be
residually finite and what is the worst possible case? The first answer that
would come to mind is when $G$ ($\neq I$) has no proper finite index subgroups
at all, but we have been dealing with groups possessing positive first Betti
number and hence infinitely many subgroups of finite index. By noting that
elements outside the commutator subgroup $G'$ cannot be in $S$, we obtain
our condition.
\begin{defn}
We say that the finitely generated group $G$ is {\bf residually useless} if
\[G'=\bigcap_{F\leq_f G}F\]
but $G$ is non-abelian.
\end{defn}
Note that by excluding $G$ being abelian, we have that $G$ residually finite
implies $G$ is not residually useless which is of course what we want. Also
$G$ has no proper finite index subgroups if and only if $G=G'$ and $G$ is
residually useless or $I$. The definition has many equivalent forms, some of
which are close to being mere rewordings, but the 
general idea is that a residually useless group cannot be distinguished from
its abelianisation if one only uses standard information about its finite
index subgroups.
\begin{prop}
Let $G$ be finitely generated and non-abelian with commutator subgroup $G'$,
abelianisation $\overline{G}=G/G'$ and let $S$ be the intersection of the
finite index subgroups of $G$. The following are equivalent:\\
(i) $G$ is residually useless.\\
(ii) $G/S$ is abelian.\\
(iii) $G$ has no non-abelian finite quotient.\\
(iv) If $a_n(G)$ denotes the number of finite index subgroups of $G$ having
index $n$ then $a_n(G)=a_n(\overline{G})$ for all $n$.\\
(v) For all $F\leq_f G$ we have $F'=G'$.\\
(vi) For all $F\leq_f G$ we have $F\cap G'=G'$.\\
(vii) For all $F\leq_f G$ we have $F'=F\cap G'$.\\
\end{prop}
\begin{proof}
The equivalence of (i) with (ii) is immediate on seeing that we always have
$S\leq G'$ if $G$ is finitely generated, and likewise with (iii) on dropping
down to a finite index normal subgroup. As for (iv), this is just using the
index preserving correspondence between the subgroups of $\overline{G}$ and 
the subgroups of $G$ containing $G'$.

As for the rest, we have that $F'\leq F\cap G'\leq G'$ whenever $F$ is a 
subgroup of $G$. If (i) holds for $G$ and $F$ is a finite index subgroup then
the intersection $R$ of the finite index subgroups of $F$ contains $S$ (which
is $G'$) but is inside $F'$ so $F'$ and $G'$ are equal, giving (v). This
immediately implies (vi) and (vii) so we just require that they in turn imply
(i). This is obvious for (vi) but not for (vii) until we spot \cite{lr} Lemma
3.2. If (i) fails then take $F\leq_f G$ and $g\in G'$ but $g\notin F$.
Dropping down to $N\leq F$ with
$N\unlhd_fG$, we have $H=N\langle g\rangle\leq_f G$ and
$H/N\cong\langle g\rangle/(N\cap\langle g\rangle)$. Thus $g\notin H'$ because
by being outside $N$ it survives under a homomorphism from $H$ to an abelian
group. But $g$ is certainly in $H\cap G'$.
\end{proof}

The importance of condition (vii) holding for $G$ 
is that we fail to pick up extra abelianisation in finite covers 
$F\leq_f G$ since $F/F'$ is just $F/F\cap G'\cong FG'/G'\leq_fG/G'$. 
In particular $\beta_1(F)=\beta_1(G)$ so $G$ is not large.

\begin{ex}
\end{ex}
(i) The Thompson group $T$ is residually useless. This group has a 2-generator
2-relator presentation with abelianisation $\z\times\z$ and its commutator
subgroup $T'$ has no proper finite index subgroups as $T'$ is infinite and
simple. Thus for $F\leq_fT$ we have $F\cap T'\leq_fT'$ giving $T'\leq F$.\\
(ii) As already mentioned, if $G$ is infinite but has no proper finite index
subgroups then $G$ is residually useless. Moreover for any such $G$ and any
finitely generated abelian group $A$ we have $\Gamma=G*A$ is residually
useless because if $N\unlhd_f\Gamma$ then $N\cap G\unlhd_f G$ so $G\leq N$.
This implies that the normal closure $C$ of $G$ is in $N$ so $\Gamma/N$ must
be abelian as it is a finite quotient of $\Gamma/C\cong A$. This also works
if $A$ is residually useless.\\
(iii) A famous example that will do for $G$ in (ii) is the Higman group $H$
with 4 generators and 4 relators, so has deficiency 0. Thus $H*H$ is residually
useless so it too has no proper finite index subgroups, given that it is 
infinite and equals its own commutator subgroup. On taking the free product of
lots of copies of $H$ and topping it off with a deficiency 1 abelian group,
namely $\z$ or $\z\times\z$, we have examples of deficiency 1 residually 
useless groups $G$ which need arbitrarily many generators and which have
$\beta_1(G)=1$ or 2.\\
(iv) We can ask if there are 1-relator groups which are residually useless.
The presentation must have 2 generators (to avoid large or cyclic cases) so
a first attempt would be the Baumslag-Solitar groups $G=B(l,m)$ with $l$ and
$m$ coprime (for if not then $G$ is large by \cite{ep2} Theorem 6). We have
in \cite{ep2} Example 3.2 a proof that $\beta_1(H)=1$ for all $H\leq_fG$ so
these are not large, but they are not residually useless either.
\begin{prop}
The Baumslag-Solitar groups are not residually useless.
\end{prop}
\begin{proof}
If $G=B(l,m)$ where by replacing generators with their inverses
we can take $l+m\geq 0$ and $l>0$ then we can see that
they surject to dihedral groups because the presentation
\[\langle x,y|xy^lx^{-1}=y^m,x^2=y^{l+m}=1,xyx^{-1}=y^{-1}\rangle\]
is a homomorphic image but the Baumslag-Solitar relation is redundant, leaving
us the dihedral group of order $2(l+m)$. This gives us a finite non-abelian
image unless $l+m=1$ or 2. In the latter case we must have $l\geq 2$ and $m\leq
-1$ as $B(1,1)$ is abelian. Now if $G$ is residually useless and $F\leq_fG$ 
then $F$ is too by Proposition 5.2 (v) so we drop to the index 2 subgroup with
exponent sum of $x$ equal to 0 modulo 2 and on rewriting we have
\[\langle t,y,z|ty^lt^{-1}=z^m,z^l=y^m\rangle\]
where $t=x^2$ and $z=xyx^{-1}$. Now as $l$ and $m$ are coprime this surjects to
\[\langle t,y,z,u|ty^lt^{-1}=z^m,u^j=1,y=u^l,z=u^m\rangle\]
where $j=l^2+m^2\geq 5$. But on adding $tut^{-1}=u^{-1}$ and $t^2=1$ to the 
presentation we see the first relation goes as before 
and we obtain our dihedral image.
\end{proof}

However in \cite{ep2} the very next example due to Edjvet and Howie is that of
an HNN extension of $B(2,3)$ given by 
\[G=\langle x,y,z|xy^2x^{-1}=y^3,zyz^{-1}=x^{-1}\rangle\]
which is a 2-generator 1-relator group shown to have the property that if
$H\leq_f G$ then $\overline{H}=\z$. This implies that $G$ is residually
useless because otherwise by Proposition 5.2 (vii)
we would have some $H$ with a homomorphism from
$\overline{H}$ to a finite index subgroup of $\overline{G}=\z$ which is
surjective but not injective. As $\z$ is Hopfian we would have $\overline{H}
\neq\z$. In fact the first example of a 1-relator residually useless group
dates back to a short paper \cite{ba} of G.\,Baumslag in 1969 entitled 
``A non-cyclic one-relator group all of whose finite quotients are cyclic'' 
with the group in question being
\[\langle a,b|a=a^{-1}b^{-1}a^{-1}bab^{-1}ab\rangle.\]

In terms of its wide application, the following is our main result on largeness
of deficiency 1 groups.
\begin{thm}
If $G$ has a deficiency 1 presentation $\langle F_n|R\rangle$ where one of the
relators is a commutator in $F_n$ then exactly one of these occurs:\\
(i) $G=\z\times\z$.\\
(ii) $G$ is residually useless with abelianisation $\z\times\z$.\\
(iii) $G$ is large.\\
In particular if $\overline{G}\neq\z\times\z$ then $G$ is large.
\end{thm}
\begin{proof} If our relator $r=uvUV$ for $u,v$ words in the generators for
$F_n$ then we can regard $r$ as the commutator of two generators simply by
adding $u$ and $v$ to the generators and their definitions to the relators,
noting that this does not change the deficiency. As we have $\beta_1(G)\geq 2$,
we have largeness by Theorem 3.2 and Proposition 3.7 unless $\overline{G}=\z
\times\z$ with $u$ and $v$ linearly independent in homology. If so then 
either we are in (i) or (ii), or by Proposition 5.2 (vii) we have $L\leq_f G$
with $\gamma$ in $L\cap G'$ but not in $L'$, 
thus granting extra homology but we need
to take care that we retain the form of our special relator $r$. We do this by
keeping track of what happens to $r$ under Reidemeister-Schreier rewriting
when we drop to a finite index subgroup. First we can assume without loss of
generality that $u$ and $v$ generate the homology as we can drop to the
appropriate finite abelian cover $H$. This works because whenever we rewrite
with respect to a subgroup of $G$, any generators of $G$ which are in the
subgroup automatically become generators of this subgroup and any relator
consisting just of these generators will survive in the presentation for $H$.
Now if $\overline{H}\neq\z\times\z$ then we are done but if not then a
failure to pick up extra abelianisation means as before that $H'=H\cap G'$
but $G'\leq H$, giving $H'=G'$. 

Let $k,l$ be the minimum positive integers such that $a=u^k$ and $b=v^l$
are in $L$. We set $N$ to be the smallest abelian cover of $H$ containing
$a$ and $b$, noting
that $N$ is an abelian cover not just of $H$ but of $G$ as well. In rewriting
for $N$ in $H$ we take a Schreier transversal of the form $u^iv^j$ (for
$0\leq i<k,0\leq j<l$). Thus $a$ and $b$ will be amongst the generators for
$N$ with our relator $uvUV$ in $H$ giving rise to $abAB$ in $N$. (It is
probably easiest to see this in two stages by first dropping to the subgroup
with exponent sum of $u$ equal to 0 mod $k$ and rewriting using the obvious
transversal $u^i$, and then doing the same with $v$.) Thus we have $G'\leq N
\unlhd_f G$ with a deficiency 1 presentation  including generators $a,b$ and
relator $abAB$. 

Finally we go from $N$ to the subgroup $L\cap N\leq_f G$ which on rewriting
will keep $a$ and $b$ because they are generators in the presentation for $N$
which also lie in $L\cap N$, and  consequently $abAB$ remains too. Now our
$\gamma\in L$ from before which is in $G'\backslash L'$ is also in $N$ as 
$G'\leq N$. But we
have a surjective but non-injective homomorphism from $L/L'$ to $L/L\cap G'
\cong\z\times\z$, thus by the Hopficity of finitely 
generated abelian groups $L$
surjects to $\z\times\z\times\z_m$ with $\gamma$ mapping onto the $\z_m$
factor. But then we can restrict this surjection to the finite index subgroup
$L\cap N$ which also
contains $\gamma$ so $L\cap N$ has the right presentation and the right
homology to obtain largeness.
\end{proof}

We can now finish off Section 4 properly and with the minimum of fuss.
\begin{co} If $G=\langle\Gamma,t\rangle$ is the mapping torus of an
endomorphism $\theta$ of the free group $\Gamma$ of rank $n$ and $\z\times\z
\leq G$ then $G=\langle x,y|xyx^{-1}=x^{\pm 1}\rangle$ or is large.
\end{co}
\begin{proof}
By the comment at the beginning of Section 4, we can assume that $\theta$ is 
really an injective endomorphism of a free group $F_m$ with $m\leq n$, putting
us in the case of Corollary 4.6 which allows us to apply Theorem 5.5 to 
$H\leq_f G$. As $\z\times\z\leq G$, we do not have $m=0$ and only the two
groups above for $m=1$. Otherwise $G$ and hence $H$ contain a
non-abelian free group for $m\geq 2$
so $H$ is not in case (i) of Theorem 5.5. By Proposition
4.3 (iii) $H$ is an injective mapping torus and so the result \cite{bs} of
Borisov and Sapir tells us that $H$ is residually finite, not residually
useless. Thus $H$ and $G$ are large.
\end{proof}

We also have a strong result for largeness of 1-relator groups. Obviously if
a 1-relator presentation has three or more generators then we have largeness, 
so we want 2-generator 1-relator large groups.
\begin{co}
If $G=\langle a,b|uvUV\rangle$ where $u$ and $v$ are any elements of $F_2=
\langle a,b\rangle$ with
$uvUV$ not equal to $abAB$, $baBA$ or their cyclic conjugates when reduced and
cyclically reduced then $G$ is large or is residually useless.
\end{co}
\begin{proof}
It is well known that $G=\z\times\z$ if and only if the relator is of the 
above form
(equivalently if and only if $u,v$ form a free basis for $F_2$) so otherwise
we are in Theorem 5.5 case (ii) or (iii).
\end{proof}

This begs the question: are there groups of the above form which are residually
useless? If not then we have a definitive result on largeness for 1-relator 
groups $G$ whose relator is a commutator; moreover on being given any such 
relator it is clear from Corollary 5.7 how one would immediately
work out whether the group is large or is $\z\times\z$. Note that $G$ has a
deficiency 1 presentation with $\beta_1(G)\geq 2$ and the only examples we know
of residually useless groups with these properties are ones of the form
$H*(\z\times\z)$ where $H\neq I$ has deficiency zero and no proper finite index
subgroups, as in Example 5.3 (iii), but these will need more than two 
generators. In fact all known examples of 2-generator 1-relator groups with
the relator a commutator are even residually finite; this is Problem (OR8) in
the problem list at \cite{nypb}. It seems likely 
that they are all not residually useless and hence large (or $\z\times\z$). 
In fact it is easy to prove this for most cases.
\begin{prop}
If $G=\langle F_2|uvUV\rangle$ then $G$ can only be residually useless if 
$u,v\notin F_2'$ with the images of $u$ and $v$ linearly independent in the
abelianisation $\z\times\z$ of $F_2$ and such that $u$ is a free basis element 
for $F_2$ or $G_u=\langle F_2|u\rangle$ is residually useless, along with the
same condition for $v$.
\end{prop}
\begin{proof}
If the images of $u$ and $v$ do not generate the homology of $F_2$ up to finite
index then we are done by Theorem 3.2. Now suppose that $G_u$ is not residually
useless or $\z$ (the latter happening if and only if $u$ is an element of a
free basis for $F_2$) then as $G$ surjects to $G_u$ we see that a non-abelian 
finite image of $G_u$ is also an image of $G$. Now swap $u$ and $v$.
\end{proof}

There are many other cases for which we can conclude that $G=\langle F_2|uvUV
\rangle$ is not residually useless and hence large. The powerful algorithm
of K.\,S.\,Brown in \cite{bks} to determine whether a 2-generator 1-relator
group is a mapping torus of an injective endomorphism of a finitely
generated group (which must necessarily be free), along with \cite{bs} proving
that such groups are residually finite, suggests that on being given a 
particular $G$ we proceed as follows: first check that $uvUV$ is not
conjugate to $abAB$ or its inverse (else $G=\z\times\z$) and that $u$ and $v$
are independent in homology (else $G$ is large). Now see if $G_u$ or $G_v$
can be shown to be not residually useless and not abelian by checking the
abelianisation of low index subgroups of each of them. If either succeeds then
$G$ is large. If not then draw out the words $u$, $v$ and $uvUV$ whilst using
Brown's algorithm to determine the BNS invariants of $G_u$, $G_v$ and $G$. If
any of these are non-empty then we have largeness of $G$ (although if the
BNS invariant of $G_u$ is non-empty then we require $G_u\neq\z$ but this is
easily checked given that $\beta_1(G_u)=1$, because if the BNS invariant
consists of both points of $S^0$ and the Alexander polynomial is 1 then
$G_u=\z$, with the same for $G_v$). If all these fail then compute the 
abelianisation of low index subgroups of $G$ and look for one which is not
$\z\times\z$ in order to obtain largeness. Needless to say, we know of no $G$
which fails all of these tests. 

We finish with a case where we can use Theorem 5.5 in practice. In \cite{ephd}
the problem of when 
\[G=\langle a,b|a^{k_1}b^{l_1}a^{k_2}b^{l_2}a^{k_3}b^{l_3}\rangle\]
for $k_1k_2k_3,l_1l_2l_3\neq 0$ and $k_1+k_2+k_3=l_1+l_2+l_3=0$ is large,
namely the case where the relator is in the derived subgroup and has free
product length 6. It was shown that $G$ is large apart from possibly the
cases of one relator and two infinite families of relators which were
\begin{eqnarray*}A^3T^2aTa^2t^3,&&\\
a^{k_1}TaTa^{k_3}t^2&\mbox{ for }&k_1+k_3=-1,\\
a^{k_1}Ta^2Ta^{k_3}t^2&\mbox{ for }&k_1+k_3=-2.
\end{eqnarray*}
Note that by \cite{ls} Proposition I.8.4 it is straightforward to tell 
whether an element $r$ in
a free group is a commutator as we must have $r=uvwUVW$ for some reduced words
$u,v,w$ where this expression involves no cancellation. Thus we can chop the
word in half and check the finite number of possibilities for $u,v,w$.
\begin{co}
The groups above are all large.
\end{co}
\begin{proof}
Writing $[u,v]$ for $uvUV$, the second case is $[a^{k_1}Ta,T^2a]$ and the third
is $[a^{k_1}Ta^2,T^2a^2]$, so Theorem 5.5 applies and we are not in case (i)
as the words are cyclically reduced. By drawing out the relators on a 2 
dimensional grid using Brown's algorithm, we conclude that the groups are
free-by-$\z$ and hence large. As for the first relator, this is not a
commutator so we turn to the computer. We find an index 4 subgroup with
abelianisation $\z^5$ and on rewriting we are told that it has a presentation
with the five generators $x=at$, $y=ta$, $z=t^2A^2$ and $a^3T$, $a^2tA$ along
with four relators including one which is $ZyxzYX$, a commutator.
\end{proof}

\end{document}